\documentclass[a4paper,12pt]{article}
\usepackage{amsmath}
\usepackage{amsthm}
\usepackage{graphicx}
\usepackage{verbatim}
\usepackage{epsfig}
\usepackage{hyperref}
\usepackage{float}
\usepackage{color}
\usepackage{fullpage}
\usepackage{amsfonts}
\usepackage{amscd}
\usepackage{amssymb}
\usepackage{graphics}
\usepackage{amsmath}
\usepackage[english]{babel}
\usepackage{bbm}
\usepackage{yfonts}
\usepackage{mathrsfs}
\usepackage{color}
\DeclareFontFamily{OML}{rsfs}{\skewchar\font'177}
\DeclareFontShape{OML}{rsfs}{m}{n}{ <5> <6> rsfs5 <7> <8> <9>
rsfs7 <10> <10.95> <12> <14.4> <17.28> <20.74> <24.88> rsfs10 }{}
\DeclareMathAlphabet{\mathfs}{OML}{rsfs}{m}{n}

\newcommand{\BE}{{\mathbb{E}}}

\newcommand{\BH}{{\mathbb{H}}}

\newcommand{\BP}{{\mathbb{P}}}

\newcommand{\BR}{{\mathbb{R}}}
\newcommand{\BS}{{\mathbb{S}}}

\newcommand{\BZ}{{\mathbb{Z}}}

\newcommand{\CB}{{\mathcal{B}}}
\newcommand{\CC}{{\mathcal{C}}}

\newcommand{\CL}{{\mathcal{L}}}

\newcommand{\CV}{{\mathcal{V}}}

\newcommand{\CX}{{\mathcal{X}}}

\newcommand{\bae}{\begin{equation}\begin{aligned}}
\newcommand{\eae}{\end{aligned}\end{equation}}

\newcommand{\cc}{\mathfrak{c}}

\newtheorem{thm}{Theorem}[section]
\newtheorem{prop}[thm]{Proposition}
\newtheorem{lem}[thm]{Lemma}

\begin{document}
\numberwithin{equation}{section} \numberwithin{figure}{section}
\title{Poisson cylinders in hyperbolic space}
\author{Erik I. Broman\footnote{Department of Mathematics, Uppsala University, Sweden. Research supported by the Swedish Research Council. E-mail: broman@math.uu.se}
\ and Johan Tykesson\footnote{Department of Mathematics, Chalmers University of Technology and Gothenburg University, Sweden. Research supported by the Knut and 
Alice Wallenberg foundation E-mail: johant@chalmers.se}}
\date{ April 5, 2015}
\maketitle
\begin{abstract}
We consider the Poisson cylinder model in $d$-dimensional hyperbolic space. 
We show that 
in contrast to the Euclidean case, there is a phase transition in 
the connectivity of the 
collection of cylinders as the intensity parameter varies. We also show 
that for any non-trivial
intensity, the diameter of the collection of cylinders is infinite.

\end{abstract}
\section{Introduction}
In the recent paper \cite{BroTyk}, the authors considered the so-called 
Poisson cylinder model in 
Euclidean space. Informally, this model can be described as a Poisson 
process $\omega$ on the space of bi-infinite lines 
in $\BR^d.$ The intensity of this Poisson process is $u$ times a normalized 
Haar measure on this space of lines.
One then places a cylinder $\cc$ of radius one around every line $L\in \omega,$
and with a slight abuse of notation, we say that $\cc \in \omega.$ 
The main result of \cite{BroTyk} was that for any $0<u<\infty$ and any 
two cylinders $\cc_1,\cc_2\in \omega,$
there exists a sequence $\cc^1,\ldots\cc^{d-2}\in \omega$ 
such that $\cc_1\cap \cc^1 \neq \emptyset, $
$\cc^1\cap \cc^2 \neq \emptyset,\ldots \cc^{d-2}\cap \cc_2 \neq \emptyset.$ 
In words, any two cylinders
in the process is connected via a sequence of at most $d-2$ other 
cylinders. Furthermore, it was proven that with probability 
one, there exists a pair of cylinders not connected in $d-3$ steps.
The result holds for any $0<u<\infty$, and therefore there is 
no connectivity phase transition.

This is in sharp contrast to what happens for other percolation models.
For example, ordinary discrete percolation (see \cite{Grimmett}), the 
Gilbert disc model (see \cite{MR}), and the Voronoi percolation model
(see \cite{BR}) all have a connectivity phase transition. A common property
that all the above listed models exhibit is something that we 
informally refer to as a ''locality property'' and can be described as follows.
Having knowledge of the configuration in some
region $A,$ gives no, or almost no, information about the configuration 
in some other region $B,$ as long as $A$ and $B$ are well separated.
For instance, in ordinary discrete percolation 
the configurations are independent if the two regions $A,B$ are disjoint, 
while for the Gilbert disc model  with fixed disc radius $r,$ 
the regions need to be at Euclidean distance at least $2r$ in order to 
have independence. For Voronoi percolation, there is a form of 
exponentially decaying dependence, i.e. the probability that the same 
cell in a Voronoi tessellation contains both points $x$ and $y$ decays 
exponentially in the distance between $x$ and $y.$

This is however not the case when dealing with the Poisson cylinder model in 
Euclidean space. Here, the dependency is polynomially 
decaying in that 
\begin{equation} \label{eqn:nonlocal}
\BP_E[B(x,1) \leftrightarrow B(y,1)]\sim d_E(x,y)^{-(d-1)},
\end{equation}
where the index $E$ stresses that we are in the Euclidean case, and where 
$\leftrightarrow$ denotes the existence of a cylinder 
${\mathfrak c}\in \omega$ connecting
$B(x,1)$ to $B(y,1)$.
Of course, this ''non-locality'' stems from the fact that the basic objects
of our percolation model are unbounded cylinders.


In Euclidean space, the non-locality property of \eqref{eqn:nonlocal} 
and the fact that the basic percolation objects (i.e. the cylinders)
are unbounded are, at least in some sense, the same thing. However, in hyperbolic
space, the result corresponding to \eqref{eqn:nonlocal} is quite different (see
Lemma \ref{lemma:prob_balls}) in that
\begin{equation} \label{eqn:hypexpdec}
\BP_H[B(x,1) \leftrightarrow B(y,1)]\sim e^{-(d-1)d_H(x,y)},
\end{equation}
(where $\BP_H$ stresses that we are in the hyperbolic case and 
$d_H$ denotes hyperbolic distance).
Since the decay is now exponential, this is a form of locality property.
Thus, by studying this model in hyperbolic space, we can study a model 
with unbounded percolation objects, but with a locality property. This
is something that does not occur naturally in the Euclidean setting. 

Before we can present our main results, we will provide a short explanation 
of our model, see Section \ref{sec:model}
below for further details. 
Consider therefore the $d$-dimensional hyperbolic space ${\mathbb H}^d$ 
for any $d\ge 2$. 
We let $A(d,1)$ be the set of geodesics in ${\mathbb H}^d$ and let 
$\mu_{d,1}$ be the unique (up to scaling) measure on $A(d,1)$ which is 
invariant under isometries. We will sometimes 
simply refer to the geodesics of $A(d,1)$ as lines.

Let $\omega$ be a Poisson point process on $A(d,1)$ with intensity $u\mu_{d,1}$, 
where $u>0$ is our parameter.
As in the Euclidean case, given a line $L\in \omega,$ we will let 
$\cc(L)$ denote the corresponding cylinder, and abuse notation
somewhat in writing $\cc\in \omega.$ 
Let 
\[
\CC:=\bigcup_{L\in \omega}\cc(L),
\]
be the {\em occupied} set and let $\CV:=\BH^d\setminus \CC$ be the \emph{vacant} set.
Furthermore, define 
\[
u_c=u_c(d):=\inf\{u\,:\,\CC\mbox{ is a.s.\ a connected set}\}.
\]
We note that by Proposition \ref{p.01law} below, we have that 
$\BP[\CC \textrm{ is connected}]\in \{0,1\}.$ 

Our main result is the following.
\begin{thm} \label{thm:main}
For any $d\ge 2$ we have $u_c(d)\in (0,\infty)$. Furthermore, for any $u>u_c,$ 
$\CC$ is connected. 
\end{thm}
\noindent {\bf Remarks:} Theorem \ref{thm:main} indicates that even though 
the cylinders are unbounded, the exponential decay of \eqref{eqn:hypexpdec}
seems to be the important feature in determining the existence of a phase
transition. The second part of the theorem is a monotonicity property, 
proving that when $u$ is so large that $\CC$ is a connected set, then 
we cannot have that $\CC$ is again disconnected for an even larger $u.$


In \cite{TT} a result similar to Theorem~\ref{thm:main} for the 
random interlacements model on 
certain non-amenable graphs was proven. The random interlacements model 
(which was 
introduced in \cite{Szn}) is a discrete percolation model 
exhibiting long-range dependence. 
However, the dependence structure for this model is very different from 
that of the Poisson cylinder model. To see this, consider three 
points $x,y,z\in \BH^d$ 
(or $\BR^d$ in the Euclidean case). If we know that there is a 
geodesic $L\in \omega$
such that $x,y\in L,$ then this will determine whether $z\in L$. 
For a random interlacement 
process, the objects studied are essentially trajectories of bi-infinite simple 
random walks, and so knowing
that a trajectory contains the points $x,y\in \BZ^d$ will give some information 
whether the trajectory contains $z\in \BZ^d,$ but not ''full'' information. Thus, 
the dependence structure is in some sense more rigid for the cylinder process.

\medskip

Knowing that $\CC$ is connected, it is natural to consider the diameter
of $\CC$ defined as follows.
For any two cylinders $\cc_a,\cc_b\in \omega,$ let 
$\textrm{Cdist}(\cc_a,\cc_b)$ be the minimal number $k$ of cylinders 
$\cc_1,\ldots\cc_k\in \omega$ such that 
\[
\cc_a\cup \cc_b \bigcup_{i=1}^k \cc_i
\]
is a connected set. If no such set exists, we say that 
$\textrm{Cdist}(\cc_a,\cc_b)=\infty.$
We then define the diameter of $\CC$ as 
\[
\textrm{diam}(\CC)=\sup\{\textrm{Cdist}(\cc_a,\cc_b):\cc_a,\cc_b\in \omega\}.
\]
Our second main result is
\begin{thm}\label{thm:infdiam}
For any $u\in (0,\infty),$ we have that 
\[
\BP[\textrm{diam}(\CC)=\infty]=1.
\]
\end{thm}
\noindent
{\bf Remark:} Of course, the result is trivial for $u < u_c.$

\medskip

When $0<u<u_c(d),$ it is natural to ask about the number  of 
unbounded components. Our next proposition addresses this.
\begin{prop} \label{prop:infnbrinfcomp}
For any $u\in(0,u_c)$ the number of infinite connected components of $\CC$
is a.s. infinite.
\end{prop}

One of the main tools will be the following discrete time particle process.
Since we believe that it may be of some independent interest, we present it here
in the introduction, along with our main result concerning it.
In essence, it behaves like 
a branching process where every particle gives rise to an infinite number
of offspring whose types can take any positive real value.

Formally, let $\xi_0, (\xi_{k,n})_{k,n=1}^\infty$ be an i.i.d. 
collection of Poisson processes
on $\BR$ with intensity measure $ue^{\min(0,x)}dx.$ Let $\zeta^0=\{0\},$ and we
think of this as the single particle in generation 0. Then, let 
$\zeta^1=\{x\geq 0:x\in \xi_{0}\}$ be the particles of generation 1,
and let $Z_{1,1}=\min\{x\in \zeta^1\}$ and inductively for any $k\geq 2,$
let $Z_{k,1}=\min\{x\in \zeta^1:x>Z_{k-1,1}\}$. Thus 
$Z_{1,1}<Z_{2,1}<\cdots$ and $\{Z_{1,1},Z_{2,1},\ldots\}=\zeta^1.$
We think of these as the offspring of $Z_{1,0}=\{0\}.$ 
In general, if $\zeta^n$ has been defined, and
$Z_{1,n}<Z_{2,n}<\cdots$ are the points in $\zeta^n,$ we let
\begin{equation} \label{eqn:particleproc}
\zeta^{k,n+1}=\bigcup_{x\in \xi_{k,n}:x+Z_{k,n}\geq 0}\{x+Z_{k,n}\},
\end{equation}
and $\zeta^{n+1}=\bigcup_{k=1}^\infty \zeta^{k,n+1}.$ We think of $\zeta^{n+1}$
as the particles of generation $n+1,$ and $\zeta^{k,n+1}$ as the offspring
of $Z_{k,n}\in \zeta^n.$ From \eqref{eqn:particleproc}, we see that 
$\zeta^{k,n+1}\subset \BR^+.$ Furthermore, conditioned on $Z_{k,n}=x,$ 
$Z_{k,n}$ gives rise to 
new particles in generation $n+1$ according to a Poisson process with intensity 
measure $d\mu_x=I(y\geq 0)ue^{-(x-y)^+}dy$ (where $I$ is an indicator function
and $(x-y)^+=\max(0,x-y)$). We let $\zeta=(\zeta^n)_{n=1}^\infty$ denote this 
particle process. We point out that in our definition,
any enumeration of the particles of $\zeta^n$ would be as good
as our ordering $Z_{1,n}<Z_{2,n}<\cdots$, 
as long as the enumeration does not depend on ''the future'', i.e.
$(\xi_{k,n+1})_{k=1}^\infty$ or such.

Informally the above process can be described as follows. 
Thinking of a particle as a point in $\BR^+$ corresponding to 
the type of that particle, it  gives rise to new points with 
a homogeneous rate {\em forward} of the position of the point, 
but at an exponentially decaying rate {\em backward} of the position 
of the point. Of course, since any individual gives rise to an infinite number
of offspring, the process will never die out. However, it can still 
die out weakly in the sense that for any $R$ there will eventually 
be no new points of type $R$ or smaller.


For any $n,$ let
\begin{equation} \label{eqn:Xndef}
X^n_{[a,b]}=\sum_{k=1}^\infty I(Z_{k,n}\in[a,b]).
\end{equation}
Thus, $X^n_{[a,b]}$ is the number of individuals in generation $n$ of type 
between $a$ and $b.$ 
We have the following theorem 
\begin{thm} \label{thm:bp}
There exists a constant $C<\infty$ such that for $u<1/4,$ and any $R<\infty$ 
\[
\sum_{n=1}^\infty\BE[X^n_{[0,R]}]<C\frac{e^{4uR}}{1-4u}<\infty.
\]
That is, $\zeta$ dies out weakly. Furthermore, 
for any $u>1/4,$
\[
\lim_{n \to \infty} \BE[X^n_{[0,R]}]=\infty.
\]
\end{thm}
Theorem \ref{thm:bp} will be used to prove that $u_c(d)>0$
(part of Theorem \ref{thm:main}) through a coupling procedure 
informally described in the
following way (see Section \ref{subsec:icp} for the formal definition). 
Consider a deterministic cylinder $\cc_0$ passing through the origin $o\in \BH^d$
and a Poisson process of cylinders in $\BH^d$ as described above. 
Let $\cc_{1,1},\cc_{2,1},\ldots$
be the set of cylinders in this process that intersect $\cc_0.$ 
These are the first generation 
of cylinders (and correspond to $\zeta^1).$ In the next step, 
we consider independent Poisson 
processes $(\omega_{k,1})_{k=1}^\infty$ and the collection of cylinders in 
$\omega_{k,1}$ that intersect $\cc_{k,1}$ (these collections will correspond to 
$(\zeta^{k,2})_{k=1}^\infty$
and the union of them corresponds to $\zeta^2).$ 
We then proceed for future generations
in the obvious way. By a straightforward coupling of this ''independent
cylinder process'' and the original one described above
(and since in every step
we use an independent process in the entire space $\BH^d$), 
we get that the set of cylinders
connected to $\cc_0$ through this procedure, will contain the set of cylinders in 
$\CC\cup \cc_0$ connected to $\cc_0.$ With some work, the 
independent cylinder process can be compared
to the particle process as indicated. 
By Theorem \ref{thm:bp}, for $u<1/4,$ the latter dies out weakly. We will
show that this implies that the number of cylinders 
(in the independent cylinder process)
connected to $\cc_0$ and intersecting $B(o,R)$ 
will be of order at most $e^{4ucR}$ where $c<\infty.$ 
However, the number of cylinders in 
$\CC$ intersecting $B(o,R)$ must be of order $e^{(d-1)R},$ 
which of course is strictly larger than $e^{4ucR}$ for 
$u>0$ small enough. Assuming that $\CC$ is connected
then leads to a contradiction.

\medskip

We end the introduction with an outline of the rest of the paper. In
Section \ref{sec:model} we give some background on hyperbolic geometry
and define the cylinder model. In Section \ref{sec:prel}, we establish some 
preliminary results on connectivity probabilities that will be useful in later sections.
In Section \ref{sec:ubamou}, we prove that $u_c(d)<\infty$ and the monotonicity part
of Theorem \ref{thm:main}. In Section \ref{s.actbt}, we prove Theorem \ref{thm:bp},
which (as described) will be a key ingredient in proving $u_c(d)>0,$ which is done 
in Section \ref{sec:lb}. In Sections \ref{sec:infdiam} and \ref{sec:inic}
we prove Theorem \ref{thm:infdiam} and Proposition \ref{prop:infnbrinfcomp} 
respectively.

\section{The model} \label{sec:model}
In this section we will start with some preliminaries of hyperbolic space
which we will have use for later, and proceed by defining
the model.
 
\subsection{Some facts about $d$-dimensional hyperbolic space}
There are many models for $d$-dimensional hyperbolic space
(see for instance \cite{Beardon},\cite{Ratcliffe} or \cite{Santalo}). 
In this paper, we prefer to consider the so-called Poincar\'e ball model.
Therefore, we consider the unit ball 
$U_d=\{x\in {\mathbb R}^d\,:\,d_E(o,x)<1\}$ (where $d_E$ denotes 
Euclidean distance)  equipped with the hyperbolic 
metric $d_H(x,y)$ given by
\begin{equation} \label{eqn:hyp-eqdistrel2}
d_H(x,y)=\cosh^{-1}
\left(1+2\frac{d_E(x,y)^2}{(1-d_E(o,x)^2)(1-d_E(o,y)^2)}\right).
\end{equation}
We refer to $U_d$ equipped with the metric $d_H$ as the Poincar\'e ball 
model of $d$-dimensional hyperbolic space, and 
denote it by ${\mathbb H}^d$.

For future convenience, we now state two well known (see for
instance Chapter 7.12 of \cite{Beardon}) rules from hyperbolic geometry. Here, 
we consider a triangle (consisting of segments of geodesics in $\BH^d$)  
with side lengths $a,b,c$ and we let 
$\alpha,\beta,\gamma$ denote the angles opposite of the segments 
corresponding to $a,b$ and $c$ respectively.\\
{\bf Rule 1:}
\begin{equation}\label{eqn:CS1}
\cosh(c)=\cosh(a)\cosh(b)-\sinh(a)\sinh(b)\cos(\gamma)
\end{equation}
{\bf Rule 2:}
\begin{equation}\label{eqn:CS2}
\cosh(c)=\frac{\cos(\alpha)\cos(\beta)+\cos(\gamma)}{\sin(\alpha)\sin(\beta)}
\end{equation}
These rules are usually referred to as hyperbolic cosine rules.





Let $\BS^{d-1}$ denote the unit sphere in $\BR^d.$ We will identify 
$\partial \BH^d$ with $\BS^{d-1}$. 
Any point $x\in\BH^d$ is then uniquely determined by the distance 
$\rho=d_H(o,x)$ of $x$ from the origin $o,$ and a point $s\in \BS^{d-1}$
by going along the geodesic from $o$ to $s$ a distance $\rho$ from $o.$
If we let $d\nu_{d-1}$ denote the solid angle element so that 
$O_{d-1}=\int_{\BS^{d-1}} d\nu_{d-1}$
is the $(d-1)$-dimensional volume of the sphere $\BS^{d-1},$ then 
the volume measure in $\BH^d$ can be expressed in hyperbolic spherical 
coordinates  (see \cite{Santalo}, Chapter 17) 
as
\[
dv_d=\sinh^{d-1}(\rho)d\rho d\nu_{d-1}.
\]
Thus, for any $A\subset \BH^d,$ the volume $v_d(A)$ can be written as
\begin{equation} \label{eqn:volume}
v_d(A)=\int_{A}\sinh^{d-1}(\rho)d\rho d\nu_{d-1}.
\end{equation}

\subsection{The space of geodesics in $\BH^d.$}\label{s.geodesics}

Let 
$A(d,1)$ be the set of all geodesics in $\BH^d$. 
As mentioned in the introduction, a geodesic $L\in A(d,1)$ will 
sometimes be referred to as a line. Although it will have no direct relevance
to the paper, we note that it is well known (see \cite{CFKP}, section 9)
that in the Poincar\'e ball model, $A(d,1)$ consists of diameters 
and boundary orthogonal 
circular segments of the unit ball $U_d.$

For any $K\subset \BH^d,$ we let $\CL_K:=\{L\in A(d,1):L\cap K\neq \emptyset\}.$ 
If $g$ is an isometry on $\BH^d$ (i.e. $g$
is a M\"obius transform leaving $U_d$ invariant, 
see for instance \cite{Ahlfors} Chapters 2 and 3),
we define $g\CL_K:=\{gL:L\in A(d,1)\}$ (where of course $gL=\{gx:x\in L\}$).
There exists a unique measure $\mu_{d,1}$ on $A(d,1)$ which is invariant 
under isometries (i.e. $\mu_{d,1}(g\CL_K)=\mu_{d,1}(\CL_K)$), 
and normalized such that $\mu_{d,1}(\CL_{B(o,1)})=O_{d-1}$ 
(see \cite{Santalo} Chapter 17 or \cite{BJST} Section 6).

For any $L\in A(d,1)$ we let $a=a(L)$ be the point on 
$L$ minimizing the distance to the origin, and define 
$\rho=\rho(L)=d_H(o,a).$ Note that $\rho=d_H(o,L).$ 
Let $\CL^+_K:=\{L\in \CL_K:a(L)\in K\}.$
According to (17.52) of \cite{Santalo}, we have that
\begin{eqnarray}\label{eqn:linevolume}
\lefteqn{\mu_{d,1}(\CL_{B(o,r)})=\mu_{d,1}(\CL^+_{B(o,r)})}\\ & &
=\frac{(d-1)O_{d-1}}{\sinh^{d-1}(1)}\int_0^r \cosh(\rho)\sinh^{d-2}(\rho)d\rho=\frac{O_{d-1}}{\sinh^{d-1}(1)}\sinh^{d-1}(r). \nonumber
\end{eqnarray}


\subsection{The process}

We consider the following space of point measures on $A(d,1)$:
\[
\Omega=\{\omega=\sum_{i=0}^{\infty}\delta_{L_i}\text{ where $L_i\in A(d,1)$, and }
\omega(\CL_A)<\infty\text{ for all compact }A\subset \BH^d\}.
\]
Here, $\delta_L$ of course denotes Dirac's point measure at $L.$

We will often use the following standard abuse of notation: 
if $\omega$ is some point measure, then we will write 
$"L\in \omega"$ instead of $"L\in \text{supp}(\omega)"$. 
We will draw an element $\omega$ from $\Omega$ according to a Poisson 
point process with intensity
measure $u\mu_{d,1}$ where $u>0$.  We call $\omega$ a \emph{(homogeneous) 
Poisson line process} of
intensity $u$ in $\BH^d.$

If $L\in A(d,1)$, we denote by ${\mathfrak c}(L,s)$ the cylinder of base 
radius $s$ centered around $L$, i.e.
\[
{\mathfrak c}(L,s)=\{x\in \BH^d\,:\,d_H(x,L)\le s\}.
\]
If $s=1$ we will simplify the notation and write ${\mathfrak c}(L,1)={\mathfrak c}(L)$. When convenient, we will write $\cc \in \omega$ instead of $\cc(L)$
where $L\in \omega.$
Recall that the union of all cylinders is denoted by ${\mathcal C}$,
\[
{\mathcal C}={\mathcal C}(\omega)=\bigcup_{L\in \omega} {\mathfrak c}(L),
\]
and that the vacant set $\CV$ is the complement $\BH^d\setminus \CC.$ 
For an isometry $g$ on ${\mathbb H}^d$ and an event $B\subset \Omega,$
we define 
$g B :=\{\omega'\in \Omega:\omega'=g\omega \textrm{ for some } \omega \in B\}.$
We say that an event $B\subset \Omega$, is invariant under isometries if 
$gB=B$ for every isometry $g.$ Furthermore, we have the following $0-1$ law. 
\begin{prop}\label{p.01law}
Suppose that $B$ is invariant under isometries. Then ${\mathbb P}[B]\in \{0,1\}$.
\end{prop}
The proof of Proposition~\ref{p.01law}  is fairly standard, so we only give a sketch
based on the proofs of Lemma $3.3$ of \cite{TykWinPre} and Lemma 2.6 of \cite{HJ}. 
Below, $\omega_{B(x,k)}$ denotes the restriction of $\omega$ to $\CL_{B(x,k)}.$ \\

\noindent
{\bf Sketch of proof.}
Let $\{z_k\}_{k\ge 1} \subset {\mathbb H}^d$ be such that for every $k\geq 1,$ $d_H(o,z_k)=e^k$,
and let $g_k$ be an isometry mapping $o$ to $z_k$. Define 
$I_{x,k}=I(\omega\in \{\BP[B|\omega_{B(x,k)}]>1/2\})$, and note that by 
L\'evy's 0-1 law, 
\[
\lim_{k\to \infty} I_{o,k}=I_B\mbox{ a.s. }
\]
Using that $B$ is invariant under isometries, it is straightforward to prove that 
the laws of $(I_B,I_{o,k})$ and $(I_B,I_{g_k(o),k})$ are  the same,
and so $I_{g_k(o),k}$ converges in probability to $I_B$. Thus,
\begin{equation}\label{e.levy2}
\lim_{k\to\infty}{\mathbb P}[I_{o,k}=I_{g_k(o),k}=1_B]=1.
\end{equation}

The next step is to prove that $I_{o,k}$ and $I_{g_k(o),k}$ are asymptotically independent, i.e.
\begin{equation}\label{e.asymptind}
\lim_{k\to \infty}|{\mathbb P}[I_{o,k}=1,I_{g_k(o),k}=1]
-{\mathbb P}[I_{o,k}=1]{\mathbb P}[I_{g_k(o),k}=1]|=0.
\end{equation}
Essentially, ~\eqref{e.asymptind} follows from the fact that when $k$ is large, 
the probability that there is any cylinder in $\omega$ which intersects both 
$B(o,k)$ and $B(g_k(o),k)$ is very small (this is why we choose $d_H(o,z_k)$ to 
grow rapidly). For this, one uses the estimate of the measure of 
lines intersecting 
two distant balls, see Lemma~\ref{lemma:measure_balls} below.

Since $I_{o,k}$ and $I_{g_k(o),k}$ are asymptotically independent, we get
\begin{equation}\label{e.levy3}
\lim_{k\to\infty} {\mathbb P}[I_{o,k}=1,I_{g_k(o),k}=0]={\mathbb P}[B](1-{\mathbb P}[B]).
\end{equation}
The only way both~\eqref{e.levy2} and~\eqref{e.levy3} can hold 
is if ${\mathbb P}[B]\in\{0,1\}$.
\fbox{}\\

We note that the laws of the random objects $\omega$, ${\mathcal C}$ and ${\mathcal V}$ 
are all invariant under isometries of ${\mathbb H}^d$.

\section{Connectivity probability estimates.} \label{sec:prel}
The purpose of this section is to establish some preliminary estimates
on connectivity probabilities, and in particular to establish 
\eqref{eqn:hypexpdec}. This result will then be used many times in the following
sections. 

For any two sets $A,B\subset \BH^d,$ we let $A \leftrightarrow B$ denote 
the event that there exists a cylinder $\cc\in \omega$ such that 
$A\cap \cc \neq \emptyset$ and 
$B\cap \cc \neq \emptyset.$ 
We have the following key estimate.
\begin{lem}\label{lemma:prob_balls}
Let $s\in(0,\infty)$. There exists two constants $0<c(s)<C(s)<\infty$ such that for any 
$x,y\in \BH^d,$ and $u\leq 2/\mu_{d,1}(\CL_{B(o,s+1)})$ we have that
\[
c(s)\,u\,e^{-(d-1)d_H(x,y)}\leq \BP[B(x,s)\leftrightarrow B(y,s)]
\leq C(s)\,u\,e^{-(d-1)d_H(x,y)}.
\]
\end{lem}
\noindent
Lemma~\ref{lemma:prob_balls} will follow easily from 
Lemmas~\ref{l.Alemma} and~\ref{lemma:measure_balls} below, and we 
defer the proof of Lemma~\ref{lemma:prob_balls} till later.

Recall that we identify ${\mathbb S}^{d-1}$ with
$\partial {\mathbb H}^d$ in the Poincar\'e ball model. Fix a half-line 
$L_{1/2}$ emanating from the origin. For $0<\theta<\pi,$ 
let ${\mathcal L}_{L_{1/2},\theta}$ be the set of 
all half-lines $L'_{1/2}$ such that $L'_{1/2}$ emanates from the origin and 
such that 
the angle between $L_{1/2}$ and $L'_{1/2}$ is at most $\theta$. 
Let $S_{\theta}(L_{1/2})$ be 
the set of all points $s\in \partial {\mathbb H}^{d}$ such that $s$ is the 
limit point of  some half-line in ${\mathcal L}_{L_{1/2},\theta}$. Then 
$S_{\theta}(L_{1/2})$ is the intersection of $\partial {\mathbb H}^{d}$ with 
a hyperspherical cap of Euclidean height $h=h(\theta)$, where 
\begin{equation}\label{e.htheta}
h(\theta)=1-\cos(\theta).
\end{equation}
The $(d-1)$-dimensional Euclidean volume of $S_{\theta}$ is given by
\begin{equation}\label{e.hyparea}
A(\theta)=\frac{O_{d-1}}{2} I_{2h-h^2}\left(\frac{d-1}{2},\frac{1}{2}\right)
\end{equation}
where $O_{d-1}$ (as above) is the  $(d-1)$-dimensional Euclidean volume of 
${\mathbb S}^{d-1}$, 
and $I_{2h-h^2}$ is a regularized incomplete beta function
(this follows from \cite{Li}, equation (1), by noting that 
$\sin^2(\theta)=2h-h^2$).

\begin{lem}\label{l.Alemma}
There are constants $0<c<C<\infty$ such that for any $\theta\le 1/10$, we have
\[
c \,\theta^{d-1} \le A(\theta)\le C\, \theta^{d-1}.
\]
\end{lem}

{\bf Proof.}
First observe that if $0\le \theta\le 1/10$, then 
$1-\theta^2/2\le \cos(\theta)\le 1-\theta^2/4$. Therefore, 
from~\eqref{e.htheta} we have
\begin{equation}\label{e.hthetaineq}
\frac{\theta^2}{4}\le h\le \frac{\theta^2}{2}
\le \frac{1}{200}\mbox{ whenever }\theta\in[0,1/10].
\end{equation}
We have
\begin{equation}\label{e.betafunc}
I_{2h-h^2}\left(\frac{d-1}{2},\frac{1}{2}\right)
=\frac{\int_{0}^{2h-h^2}t^{(d-1)/2-1}(1-t)^{1/2-1}\,dt}
{\int_{0}^{1}t^{(d-1)/2-1}(1-t)^{1/2-1}\,dt}.
\end{equation}
The denominator in~\eqref{e.betafunc} is a dimension-dependent constant. 
Furthermore, if $0\le h \le 1/8$, then $0\le 2h-h^2\le 1/4$, and if
$0\le t \le 1/4$, then $1\le 1/\sqrt{1-t}\le 2$. Hence, for $h\le 1/8$,
$$C_1 \,\int_{0}^{2h-h^2}t^{(d-1)/2-1}\,dt 
\le I_{2h-h^2}\left(\frac{d-1}{2},\frac{1}{2}\right)
\le C_2 \,\int_{0}^{2h-h^2}t^{(d-1)/2-1}\,dt$$
which after integration gives
$$ C_3(2h-h^2)^{\frac{d-1}{2}} \le I_{2h-h^2}\left(\frac{d-1}{2},\frac{1}{2}\right)\le C_4 (2h-h^2)^{\frac{d-1}{2}}.$$
Hence, for $h\le 1/8$,
\begin{equation} \label{e.hineq}
C_3 h^{(d-1)/2}\le I_{2h-h^2}\left(\frac{d-1}{2},\frac{1}{2}\right)\le C_4 (2 h)^{(d-1)/2}.
\end{equation} 
The lemma now follows from~\eqref{e.hyparea},~\eqref{e.hthetaineq} and~\eqref{e.hineq}. \fbox{}\\

\begin{lem} \label{lemma:measure_balls}
Let $s\in(0,\infty)$. There exists two constants $0<c(s)<C(s)<\infty$ such that for any $x,y\in \BH^d,$ 
we have that
\[
c(s)\, e^{-(d-1)d_H(x,y)}\leq \mu_{d,1}(\CL_{B(x,s)}\cap \CL_{B(y,s)})
\leq C(s)\, e^{-(d-1) d_H(x,y)}.
\]
\end{lem}
{\bf Proof.} 
For convenience, we perform the proof in the case $s=1$. The general case is 
dealt with in the same way. The proof is somewhat similar to the proof of 
Lemma $3.1$ in \cite{TykWinPre}. Recall that we use the Poincar\'e ball 
model, and keep in 
mind that $\partial {\mathbb H}^d$ is identified with ${\mathbb S}^{d-1}$. 
Let $R=d_H(x,y)$ and without loss of generality assume that $x=o$ and 
so $y\in\partial B(o,R).$ We can assume that $R>2$ as the case $R\le 2$ follows by 
adjusting the constants $c,C$.
For any $R\in(0,\infty]$ and $A\subset \partial B(o,R)$, let 
\[
\tau_R(A):=\mu_{d,1}(\CL_{B(o,1)}\cap \CL_{A}).
\]  
The projection $\Pi_{\partial {\mathbb H}^{d}}(A)$ of $A$ onto 
$\partial {\mathbb H}^{d}$ is defined as the set of all 
points $y$ in $\partial {\mathbb H}^d$ for which there is a half-line emanating from $o$, 
passing through $A$ and with its end-point at infinity at $y$.

We now argue that 
\begin{equation}\label{e.rotinv}
\mu_{d,1}(\CL_{B(o,1)})\sigma_R(A)\le\tau_R(A)\le 2\mu_{d,1}(\CL_{B(o,1)})\sigma_R(A),
\end{equation}
where $\sigma_R$ is the unique rotationally invariant probability measure 
on $\partial B(o,R)$. Here, $\sigma_{\infty}$ is the rotationally 
invariant probability measure on $\partial {\mathbb H}^d$, which is just 
a constant multiple of the Lebesgue measure on ${\mathbb S}^{d-1}$. 
For $A\subset \partial B(o,R)$, let $N_A(\omega)$ denote the number of 
points in $A$ that are intersected by some line in 
${\mathcal L}_{B(o,1)}\cap \omega$.  If $L\in \CL_{B(o,1)}\cap \CL_{A}$ 
then $L$ intersects $A$ at one or two points. Hence 
\begin{equation}\label{e.Neq}
N_A(\omega)/2\le \omega({\mathcal L}_{B(o,1)}\cap {\mathcal L}_A)\le N_A(\omega).
\end{equation}
In addition, every line intersecting $B(o,1)$ intersects $\partial B(o,R)$ 
exactly twice. Hence, 
\begin{equation}\label{e.Neq2}
N_{\partial B(o,R)}(\omega)=2\omega({\mathcal L}_{B(o,1)}).
\end{equation}
For $A\subset \partial B(o,R)$ define $\rho_R(A)={\mathbb E}[N_A(\omega)]$. 
Taking expectations in~\eqref{e.Neq} we obtain 
\begin{equation}\label{e.rhoeq}
\rho_R(A)/2\le u\,\tau_R(A) \le \rho_R(A).
\end{equation}

It is easily verified that $\rho_R(A)$ is invariant under rotations. 
Hence, $\rho_R$ is a constant multiple of $\sigma_R$. Taking expectations 
in~\eqref{e.Neq2}, we obtain 
\[
\rho_R(\partial B(o,R))=  2 u  \mu_{d,1}(\CL_{B(o,1)}),
\]
from which it follows that 
\begin{equation}\label{e.rhoeq2}
\rho_R(\cdot)=2 u \mu_{d,1}(\CL_{B(o,1)})\sigma_R(\cdot).
\end{equation}
Combining~\eqref{e.rhoeq} and~\eqref{e.rhoeq2} we obtain~\eqref{e.rotinv}. 
Since $\sigma_R(A)=\sigma_{\infty}(\Pi_{\partial {\mathbb H}^{d}}(A))$, this gives

\begin{equation} \label{eqn:A'}
\mu_{d,1}(\CL_{B(o,1)})\sigma_{\infty}(\Pi_{\partial {\mathbb H}^{d}}(A))
\le\tau_R(A)\le 2\mu_{d,1}(\CL_{B(o,1)})\sigma_{\infty}(\Pi_{\partial {\mathbb H}^{d}}(A)).
\end{equation}

Having proved \eqref{e.rotinv} and \eqref{eqn:A'}, we
now proceed to prove the lower bound. We observe that 
\begin{equation*}
\CL_{B(o,1)}\cap \CL_{B(y,1)}\supset  \CL_{B(o,1)}\cap \CL_{B(y,1)\cap \partial B(o,R)}.
\end{equation*}
Hence, in view of~\eqref{e.rotinv}, we need to estimate $\sigma_R(E)$ from below,
where $E=B(y,1)\cap \partial B(o,R)$. 
Let $L_1$ be any line containing $o$ 
and intersecting $\partial B(y,1)\cap\partial B(o,R)$, and let $L_y$ be
the line intersecting $o$ and $y.$
  Denote the angle 
between $L_1$ and $L_y$ by $\theta=\theta(R)$.  
Observe that $\Pi_{\partial {\mathbb H}^{d}}(E)$ is the intersection 
of $\partial {\mathbb H}^{d}$ and 
a hyperspherical cap of Euclidean height $1-\cos(\theta)$, and so we need
to find bounds on $\theta$.

Applying \eqref{eqn:CS1} to the triangle defined by $L_1\cap B(o,R)$, 
the line segment between $o$ and $y$, and the line segment between 
$L_1\cap \partial B(o,R)$ and $y$, we have
\begin{equation}\label{e.thetaeq1}
\cosh(1)=\cosh^2(R)-\sinh^2(R)\cos(\theta).
\end{equation}
Solving~\eqref{e.thetaeq1} for $\theta$ gives
\begin{equation*}
\theta=\arccos\left(1-\left(\frac{\cosh(1)-1}{\sinh^2(R)}\right)\right).
\end{equation*}
Observe that for any $0\leq x\leq 1,$
\[
\arccos(1-x)=\arcsin(\sqrt{2x-x^2})\ge \arcsin(\sqrt{x})\ge  \sqrt{x}.
\]
Hence for $R\ge 1$,
\begin{equation} \label{eqn:thetaR}
\theta\geq \frac{C}{\sinh(R)}\geq C e^{-R}.
\end{equation}
By Lemma~\ref{l.Alemma} , we have
\begin{equation}\label{eqn:E'}
\sigma_\infty(\Pi_{\partial {\mathbb H}^{d}}(E))\ge c\, \theta^{d-1},
\end{equation}
and so the lower bound follows by combining \eqref{eqn:A'}, \eqref{eqn:thetaR}
and \eqref{eqn:E'}.

We turn to the upper bound.
Let $y'$ be the point on $\partial B(y,1)$ closest to the origin, and 
let $H$ be the
$(d-1)$-dimensional hyperbolic space orthogonal to $L_y$ and containing $y'.$
Let $\Pi_{\partial {\mathbb H}^{d}}(H)
\subset \partial \BH^d$ be the projection of $H$ onto 
$\partial \BH^d$. 
Since for any $z\in H,$ $d_H(y,z)\geq d_H(y,y')$ we get that
\begin{equation*}
\CL_{B(o,1)}\cap \CL_{B(y,1)}\subset  
\CL_{B(o,1)}\cap \CL_{\Pi_{\partial {\mathbb H}^{d}}(H)}.
\end{equation*}
Next we find an upper bound of $\sigma_{\infty}(\Pi_{\partial {\mathbb H}^{d}}(H))$, 
which will imply the 
upper bound of $\mu_{d,1}(\CL_{B(o,1)}\cap \CL_{B(y,1)})$. Let $L_2$ be 
any geodesic in $H,$ and let $s$ and $s'$ 
be the two end-points at infinity of $L_2$. Let $L_3$ be the half-line 
between $0$ and $s$, and let $\gamma=\gamma(R)$ be the angle between 
$L_3$ and $L_y$.  Applying \eqref{eqn:CS2} 
to the triangle defined by $L_3$, the half-line between $s$ and $y'$ and 
the line-segment between $0$ and $y'$, we obtain
$$\cos(0)=-\cos(\pi/2)\cos(\gamma)+\sin(\pi/2)\sin(\gamma)\cosh(R-1),$$
which gives
$$1=\sin(\gamma)\cosh(R-1).$$
Observe that we here applied~\eqref{eqn:CS2} to an infinite triangle, 
which can be justified by a limit argument.
Hence
\[
\gamma=\arcsin\left(\frac{1}{\cosh(R-1)}\right).
\]
Observe that $\arcsin(x)\le 2x$ for every $0\leq x\leq 1,$ so that 
\begin{equation} \label{eqn:gammaR}
\gamma\leq \frac{2}{\cosh(R-1)}\le C e^{-R}.
\end{equation}
We observe that $\Pi_{\partial {\mathbb H}^{d}}(H)$ is the intersection 
between a hyperspherical cap 
of Euclidean height $1-\cos(\gamma)$ and $\partial {\mathbb H}^d$. 
Hence, according to Lemma~\ref{l.Alemma},
\begin{equation} \label{eqn:AR}
\sigma_{\infty}(\Pi_{\partial {\mathbb H}^{d}}(H))\le c \gamma^{d-1}.
\end{equation}
The upper bound follows by combining \eqref{eqn:A'}, \eqref{eqn:gammaR}
and \eqref{eqn:AR}, which concludes the proof.
\fbox{}\\

\begin{lem}\label{l.diffballs}
Suppose $d_H(x,y)=R$ and that $r,s\in(0,\infty)$. There is a constant $c(d,s)<\infty$ such that 
if $R>r+s$, then
$$\mu_{d,1}(\CL_{B(x,s)}\cap \CL_{B(y,r)})\le c(d,s)\exp(-(d-1)(R-r)).$$
\end{lem}
{\bf Proof.}
The proof is nearly identical to the proof of the upper bound in 
Lemma~\ref{lemma:measure_balls}, and therefore we leave the details to 
the reader. 
\fbox{}\\

\noindent
We can now prove Lemma~\ref{lemma:prob_balls}.\\
\noindent
{\bf Proof of Lemma~\ref{lemma:prob_balls}.}
We perform the proof in the case $s=1$ as the general case follows similarly. 
First observe that 
\[
\{B(x,1)\leftrightarrow B(y,1)\}=\{\omega(\CL_{B(x,2)}\cap \CL_{B(y,2)})\ge 1\}.
\] 
Using that $1-e^{-x}\leq x$ for $x\geq 0,$ we have that
\begin{eqnarray*}
\lefteqn{\BP[B(x,1)\leftrightarrow B(y,1)]
=1-\BP[B(x,1) \not \leftrightarrow B(y,1)]}\\
& & =1-e^{-u\mu_{d,1}(\CL_{B(x,2)}\cap \CL_{B(y,2)})}
\leq u\mu_{d,1}(\CL_{B(x,2)}\cap \CL_{B(y,2)})
\leq Cu e^{-(d-1)d_H(x,y)}
\end{eqnarray*}
by Lemma \ref{lemma:measure_balls} with $C$ as in the same lemma. 

Using that $1-e^{-x}\geq x/2$ if 
$x\leq 2,$ and that $u\mu_{d,1}(\CL_{B(x,2)}\cap \CL_{B(y,2)})\leq u\mu_{d,1}(\CL_{B(x,2)})=u\mu_{d,1}(\CL_{B(o,2)})\leq 2$ by assumption,
we get as above that
\[
\BP[B(x,1)\leftrightarrow B(y,1)] 
\geq \frac{u\mu_{d,1}(\CL_{B(x,2)}\cap \CL_{B(y,2)})}{2}
\geq cue^{-(d-1)d_H(x,y)}
\]
by again using Lemma \ref{lemma:measure_balls} and letting $c$ be 
half of that of Lemma \ref{lemma:measure_balls}.
\fbox{}\\

\section{Proof of $u_c<\infty$ and monotonicity of uniqueness} \label{sec:ubamou}

We start by proving the monotonicity of uniqueness.
For convenience, in this section we denote by $\omega_u$ a Poisson 
line process with intensity $u$.  In addition, 
we will let ${\mathbb E}$ and ${\mathbb P}$ denote expectation and 
probability measure for several Poisson processes simultaneously.
Recall also that $A(d,1)$ is the set of all geodesics in $\BH^d.$
\begin{lem} \label{lemma:mon}
If for $u_1>0$ ${\mathbb P}[{\mathcal C}(\omega_{u_1})\mbox{ is connected}\,]=1,$ 
then ${\mathbb P}[{\mathcal C}(\omega_{u_2})\mbox{ is connected}\,]=1$ for every $u_2>u_1.$
\end{lem}
{\bf Proof.}
It is straightforward to show that for any $L\in A(d,1)$,
$\mu_{d,1}({\mathcal L}_{\mathfrak{c}(L)})=\infty.$ Hence, for any $L\in A(d,1),$
\begin{equation}\label{e.easytoshow}
{\mathbb P}[\mathfrak{c}(L)\cap {\mathcal C}(\omega_{u_1})\neq\emptyset]=1.
\end{equation} 
Let $u'=u_2-u_1$ and let $\omega_{u'}$ be a Poisson line process of intensity 
$u'$, independent of $\omega_{u_1}.$  By the Poissonian nature of the 
process, ${\mathcal C}(\omega_{u_2})$ has the same law as 
${\mathcal C}(\omega_{u_1})\cup {\mathcal C}(\omega_{u'})$. Hence it 
suffices to show that the a.s. connectedness of 
${\mathcal C}(\omega_{u_1})$ implies the a.s. connectedness of 
${\mathcal C}(\omega_{u_1})\cup {\mathcal C}(\omega_{u'})$. To 
show this, it suffices to show that a.s., every line in $\omega_{u'}$ 
intersects ${\mathcal C}(\omega_{u_1})$. To this end, for $L\in A(d,1)$, 
define the event $S(L)= \{\mathfrak{c}(L)\cap {\mathcal C}(\omega_{u_1})\neq \emptyset\}$. 
Then let
\[
D:=\cap_{L\in \omega_{u'}}S(L).
\]
We will show that ${\mathbb P}[D^c]=0$ and we start by observing that 
\[
{\mathbb P}[D^c]=\BP\left[\cup_{L\in \omega_{u'}}S(L)^c\right]
\leq {\mathbb E}\left[\sum_{L\in \omega_{u'}}I(S(L)^c )\right].
\]
For clarity, we let $\BE^{\omega_{u'}}$ and $\BE^{\omega_{u_1}}$
denote expectation with respect to the processes $\omega_{u'}$ and 
$\omega_{u_1}$ respectively, and we will let $\BE$ denote expectation 
with respect to $\omega_{u'}\cup \omega_{u_1}$. We use similar notation
for probability. We then have that, 
\begin{eqnarray*}
\lefteqn{{\mathbb E}\left[\sum_{L\in \omega_{u'}}I( S(L)^c)\right]}\\
& & ={\mathbb E}^{\omega_{u'}}
\left[ {\mathbb E}^{\omega_{u_1}}\left[\sum_{L\in \omega_{u'}}I( S(L)^c)\,
\bigg|\,\omega_{u'}\right]\right]
={\mathbb E}^{\omega_{u'}}
\left[ \sum_{L\in \omega_{u'}}{\mathbb E}^{\omega_{u_1}}\left[I( S(L)^c)\,
|\,\omega_{u'}\right]\right] \\
& & ={\mathbb E}^{\omega_{u'}}
\left[ \sum_{L\in \omega_{u'}}\BP^{\omega_{u_1}}\left[S(L)^c\,| \,\omega_{u'}\right]\right]
={\mathbb E}^{\omega_{u'}}
\left[ \sum_{L\in \omega_{u'}}\BP^{\omega_{u_1}}\left[S(L)^c\right]\right]=0,
\end{eqnarray*}
where we use the independence between $\omega_{u'}$ and $\omega_{u_1}$ in the 
penultimate equality and that $\BP\left[S(L)^c\right]=0$ which 
follows from \eqref{e.easytoshow}. This finishes the proof of the proposition.
\fbox{}\\

The aim of the rest of this section is to prove the following proposition, which
is a part of Theorem \ref{thm:main}
\begin{prop}\label{p.cylconn}
For any $d\ge 2$,
\begin{equation*}
u_c(d)<\infty.
\end{equation*}
\end{prop}

In order to prove Proposition \ref{p.cylconn}, we will need some preliminary
results and terminology. Recall the definition of  
${\mathcal L}^+_A$ for $A\subset {\mathbb H}^d$ and the definitions of 
$a(L)$ and $\rho(L),$ all from Section \ref{s.geodesics}.
Using the line process $\omega$, we define a point 
process $\tau$ in $\BH^d$ as follows:
\[
\tau=\tau(\omega):=\sum_{L\in\omega} \delta_{a(L)}.
\]
In other words, $\tau$ is the point process induced by the points that 
minimize the distance between the origin and the lines of $\omega$.
We observe that since $\omega$ is a Poisson process, it follows 
that $\tau$ is also a Poisson process (albeit inhomogeneous).
We will consider a percolation model with balls in place of cylinders,
using $\tau$ as the underlying point process.
Our aim is to prove that $\CV$ does not percolate for $u<\infty$ 
large enough by analyzing this latter model.
For this, we will need Lemma \ref{l.closepointbound}, 
which provides a uniform bound (in $z\in \BH^d$) of the probability that a 
point of $\tau$ falls in the ball of radius $1/2$ centered at $z\in \BH^d.$
Before that, we present
the following lemma, which will be useful on several occasions.

\begin{lem}\label{lemma:D}
There exists a set $D$ of points in 
${\mathbb H}^d$ with the following properties:
\begin{enumerate}
\item $d_H(z,D)\le1/2$ for all $z\in {\mathbb H}^d$.
\item If $x,y\in D$ and $x\neq y$, then $d_H(x,y)\ge 1/2$.
\end{enumerate}
Furthermore, for any such set, there exist constants $0<c_1(d)<c_2(d)<\infty$ so that 
for any $x\in \BH^d,$ and $r\geq 1,$
\begin{equation}\label{e.duniform}
c_1(d) v_d(B(o,r)) \leq |D\cap B(x,r)| \leq c_2(d) v_d(B(o,r+1)).
\end{equation}
\end{lem}
\noindent
{\bf Proof.}
We give an explicit construction of the set $D$. First let 
$D_1=\{o\}$ and $E_1=\{x\in \BH^d:d_H(o,x)=1/2\}$, and define
$D_2=D_1\cup \{x_1\}$ where $x_1$ is any point in $E_1.$
Inductively, having defined $D_n,$ we let 
$E_n=\{x\in \BH^d:d_H(D_n,x)=1/2\}$ and define 
$D_{n+1}=D_n\cup \{x_n\}$ where $x_n$ is any point in $E_n$
such that $d_H(o,x_n)=d_H(o,E_n)$ which exists by compactness
of the set $E_n.$
Finally we let $D=\cup_{n=1}^\infty D_n.$ By construction, any 
two points in $D$ will then satisfy condition 2. 
 Assume now
that there exists a point $z\in \BH^d$ such that $d_H(z,D)>1/2,$
and let $m$ be any integer such that $d_H(o,x_m)\geq d_H(o,z).$
Since $d_H(z,D_m)\geq d_H(z,D)>1/2$ we have that
\begin{equation} \label{eqn:dzDm}
d_H\left(z,\bigcup_{x\in D_m}B(x,1/2)\right)>0.
\end{equation}
Let $S_{z}$ be the line segment from $o$ to $z,$ and observe that 
since $o\in \bigcup_{x\in D_m}B(x,1/2),$ there must be 
some point $s=s(E_m,z)$ 
belonging to $S_{z}\cap E_m.$ Because of \eqref{eqn:dzDm}, we see that 
for some $\epsilon>0,$ we have that $d_H(o,z)=d_H(o,s)+\epsilon$ and so 
we get that 
\[
d_H(o,z)=d_H(o,s)+\epsilon\geq d_H(o,E_m)+\epsilon= d_H(o,x_m)+\epsilon>d_H(o,x_m),
\]
leading to a contradiction.

We now turn to \eqref{e.duniform} and start with the upper bound. 
Let $y_1,\ldots,y_N$ be an enumeration of $D\cap B(x,r)$. By construction, 
the balls $B(y_k,1/5)$ are all disjoint, and so $Nv_d(B(o,1/5))\leq v_d(B(o,r+1))$
from which the upper bound follows with $c_2=1/v_d(B(o,1/5)).$

For the lower bound, it suffices to observe that from the construction we have that 
\[
B(o,r)\subset \bigcup_{k=1}^N B(y_k,1),
\]
so that $N\geq v_d(B(o,r))/v_d(B(o,1))$. Hence, the lower bound follows with $c_1=1/v_d(B(o,1)).$
\fbox{}\\

\begin{lem}\label{l.closepointbound}
There is a constant $c(d)>0$ such that for any $z\in {\mathbb H}^d$,
$$\mu_{d,1}({\mathcal L}_{B(z,1/2)}^+)\ge c.$$
\end{lem}
\noindent
{\bf Proof.}
We first claim that there is a constant $c_1=c_1(d)\in (0,\infty)$ 
such that for any $r\ge 0$, the shell $B(o,r+1/4)\setminus B(o,(r-1/4)^+)$ can be covered 
by at most $N_r=\lceil c_1 e^{(d-1) r}\rceil$ balls of radius $1/2$ centered 
in $\partial B(o,r)$.  For this, we observe that by modifying the proof of 
Lemma \ref{lemma:D}, we can obtain a set of points $E\subset{\mathbb H}^d$ 
with the properties that $d(x,E)\le 1/4$ for all $x\in {\mathbb H}^d$ and 
$|E\cap B(o,r+1/2)|\le c \nu_d (B(o,r+3/2))$ for some constant $c<\infty$ and all $r\geq 1$.  
Let $E_r=E\cap (B(o,r+1/2)\setminus B(o,(r-1/2)^+)$. Since $d(x,E)\le 1/4$ for 
all $x\in {\mathbb H}^d$ we have 
\[
B(o,r+1/4)\setminus B(o,(r-1/4)^+)\subset \bigcup_{x\in E_r}B(x,1/4).
\]
For $x\in E_r$ let $x'$ be the point on $\partial B(o,r)$ minimizing the 
distance between $x$ and $\partial B(o,r),$ and let $E_r'\subset \partial B(o,r)$ denote the 
collection of all such $x'$. Since $d(x,x')\le 1/4$ we have $B(x,1/4)\subset B(x',1/2)$. Hence
\[
B(o,r+1/4)\setminus B(o,(r-1/4)^+)\subset \bigcup_{x'\in E_r'}B(x',1/2).
\]
The claim follows, since $|E_r'|\leq |E\cap B(o,r+1/2)|\le c \nu_d(B(o,r+3/2))\le c' e^{(d-1) r}$.

Now fix $z\in {\mathbb H}^d$ and let $r:=d_H(o,z)$. 
The $\mu_{d,1}$-measure of lines that have their closest point to the origin inside the 
shell $B(o,r+1/4)\setminus B(o,(r-1/4)^+)$ is given by 
\begin{eqnarray}\label{e.nreq}
\lefteqn{\mu_{d,1}({\mathcal L}_{B(o,r+1/4)}\setminus {\mathcal L}_{B(o,(r-1/4)^+)})} \\ 
& & =\mu_{d,1}({\mathcal L}_{B(o,r+1/4)})-\mu_{d,1}({\mathcal L}_{B(o,(r-1/4)^+)}) \nonumber\\ 
& &=C(d)(\sinh^{d-1}(r+1/4)-\sinh^{d-1}((r-1/4)^+)\ge C'(d)\, e^{(d-1)r},\nonumber
\end{eqnarray}
where the second equality uses \eqref{eqn:linevolume} with $C(d)=(d-1)/\sinh^{d-1}(1).$
Let $(x_i)_{i=1}^{N_r}$ be a collection of points in $\partial B(o,r)$ such that 
\begin{equation}\label{e.balleq}
B(o,r+1/4)\setminus B(o,(r-1/4)^+)\subset \cup_{i=1}^{N_r} B(x_i,1/2).
\end{equation}
From~\eqref{e.nreq} and~\eqref{e.balleq} we obtain
\begin{equation}\label{e.nreq2}
C'(d) e^{(d-1) r}\le \sum_{i=1}^{N_r}\mu_{d,1}({\mathcal L}^+_{B(x_i,1/2)})=N_r \mu_{d,1}({\mathcal L}^+_{B(z,1/2)}),
\end{equation}
where we used that $\mu_{d,1}$ is invariant under rotations in the last equality. From~\eqref{e.nreq2} we conclude that
\[
\mu_{d,1}({\mathcal L}^+_{B(z,1/2)}) \ge C'(d)e^{(d-1)r} / N_r\ge c(d)>0,
\]
finishing the proof of the lemma. \fbox{}\\


\begin{prop}\label{p.pbunique}
For any $d\ge 2$, the set ${\mathcal V}$ does not percolate if 
$u$ is large enough.
\end{prop}
\noindent
{\bf Proof.}
The proof follows the proof of Lemma $6.5$ in \cite{BS} quite closely. 
Let 
\[
{\mathcal W}:=\left(\bigcup_{x\in \tau} B(x,1)\right)^c.
\]
Then it is clear that ${\mathcal W}\supset {\mathcal V}$ so it suffices 
to show that ${\mathcal W}$ does not percolate when $u$ is large.

For $z\in {\mathbb H}^d$ let $Q(z)$ be the event that $z$ is within 
distance $1/2$ from ${\mathcal W}$. Then $Q(z)$ is determined by 
$\tau \cap B(z,3/2)$ so that $Q(z)$ and $Q(z')$ are independent if $d_H(z,z')\ge 3$.
Let $A$ be the event that $o$ belongs to an infinite component of 
${\mathcal W}$.  If $A$ occurs, then there exists an infinite 
continuous curve $\gamma\,:\,[0,\infty)\to{\mathcal W}$ with the 
properties that $\gamma(0)=o$ and $d_H(o,\gamma(t))\to\infty$ as 
$t\to\infty$. Let $t_0=0$ and $y_0=o$, and for $k\ge 1$ let 
inductively $t_k=\sup\{t\,:\,d_H(\gamma(t),y_{k-1})=6\}$ and $y_k=\gamma(t_k)$. 
For each $k$, let $y_k'$ be a point in $D$ which minimizes the distance to 
$y_k$. By definition $d_H(y_j,y_k)\ge 6$ if $j\neq k$, and since
$d_H(y_j,y_j')\le 1/2$ and $d_H(y_k,y_k')\le 1/2$, we get $d_H(y_j',y_k')\ge 5$ if 
$j\neq k$. Since $d_H(y_k,y_{k+1})=6$ we also have $d_H(y_k',y_{k+1}')\le 7$. 
Observe that since $y_k\in {\mathcal W}$, the event $Q(y_k')$ occurs. 

Let $D$ be as in Lemma \ref{lemma:D}, and
let $\CX_n$ be the set of sequences $x_0,...,x_n$ of points in $D$ such 
that $d_H(o,x_0)\le 1/2$, $d_H(x_n,x_{n+1})\le 7$ and 
$d_H(x_j,x_k)\ge 5$ if $j\neq k$. Furthermore, let $N_n$ denote the number of such 
sequences. We have that 
\begin{equation}\label{e.AQeq1}
{\mathbb P}[A]\le \sum_{(x_0,...,x_n)\in \CX_n} {\mathbb P}[Q(x_0)\cap...\cap Q(x_n)],
\end{equation}
and that 
\begin{equation}\label{e.AQeq2}
N_n\le \sup\{|D\cap B(z,7)|^{n+1}\,:\,z\in 
{\mathbb H}^d\}\stackrel{~\eqref{e.duniform}}{\le}c_1(d)^{n+1}
\end{equation}
for some constant $c(d)<\infty$.
By independence, 

\begin{equation}\label{e.AQeq3}
{\mathbb P}[Q(x_0)\cap ... \cap Q(x_n)]=\Pi_{i=0}^{n}{\mathbb P}[Q(x_i)].
\end{equation}
Observe that if $\tau(B(z,1/2))\ge 1$, then 
$B(z,1/2)\subset {\mathcal W}^c$. Hence we have
\begin{eqnarray}\label{e.AQeq4}
\lefteqn{{\mathbb P}[Q(z)]
={\mathbb P}[B(z,1/2)\cap {\mathcal W}\neq\emptyset]} \\ 
& & \le{\mathbb P}[\tau(B(z,1/2))=0]
=e^{-u \mu_{d,1}({\mathcal L}^{+}_{B(z,1/2)})} \le e^{-u c(d)},\nonumber
\end{eqnarray}
where the last inequality follows from Lemma~\ref{l.closepointbound}.
From~\eqref{e.AQeq1}, ~\eqref{e.AQeq2}, ~\eqref{e.AQeq3}, and ~\eqref{e.AQeq4} it follows that
\begin{equation}\label{e.Aprobupper}
{\mathbb P}[A]\le (c_1(d) e^{-u c(d)})^{n+1}\to 0.
\end{equation}
 as $n\to \infty$ if $u<\infty$ is large enough. 
 We conclude that ${\mathbb P}[A]=0$ for $u$ large enough but finite.
\fbox{}\\

\noindent
We can now prove Proposition \ref{p.cylconn}.\\
\noindent
{\bf Proof of Proposition~\ref{p.cylconn}: }If ${\mathcal C}$ is 
disconnected, then it consists of more than one infinite connected 
component. Since any two disjoint infinite components of ${\mathcal C}$ 
must be separated by some infinite component of ${\mathcal V}$, 
we get that the disconnectedness 
of ${\mathcal C}$ implies that ${\mathcal V}$ percolates. According 
to Proposition~\eqref{p.pbunique}, there is no percolation in 
${\mathcal V}$ when $u$ is large enough. Hence ${\mathcal C}$ is connected 
when $u$ is large enough.\fbox{}\\

\section{Proof of Theorem \ref{thm:bp}.}\label{s.actbt}

Before we can prove Theorem \ref{thm:bp}, we will need 
to do some preliminary work. 
To that end, let $\{c_{k,n}\}_{n\geq 0, -1\leq k\leq n}$ be defined by letting
$c_{0,0}=c_{0,1}=c_{1,1}=1$ and $c_{-1,n}=0$ for every $n$ and then 
inductively for every $0\leq k\leq n$ letting
\begin{equation} \label{eqn:cattri}
c_{k,n}:=\sum_{l=k-1}^{n-1}c_{l,n-1},
\end{equation}
where we define $c_{n+1,n}=0.$ Note that by this definition,  
$c_{k,n}=c_{k-1,n-1}+c_{k+1,n}.$
These numbers constitute (a version) of the Catalan triangle, 
and it is easy to verify that 
\begin{equation} \label{eqn:ckn}
c_{k,n}
=\frac{(2n-k)!(k+1)}{(n-k)!(n+1)!}
=\frac{k+1}{n+1}{2n-k \choose n}
\end{equation}
for every $n$ and $0\leq k \leq n.$
This follows by using that if \eqref{eqn:ckn} holds for 
$c_{k-1,n-1}$ and $c_{k+1,n}$, we get that 
\begin{eqnarray*}
\lefteqn{c_{k,n}=c_{k-1,n-1}+c_{k+1,n}}\\
& & =\frac{(2n-k-1)!k}{(n-k)!n!}
+\frac{(2n-k-1)!(k+2)}{(n-k-1)!(n+1)!}\\
& & =\frac{(2n-k-1)!k(n+1)+(2n-k-1)!(k+2)(n-k)}{(n-k)!(n+1)!}\\
& & =\frac{(2n-k-1)!(2kn-k+2n-k^2)}{(n-k)!(n+1)!}
=\frac{(2n-k)!(k+1)}{(n-k)!(n+1)!}.
\end{eqnarray*}
By an induction argument, we see that 
\eqref{eqn:ckn} holds for every $0\leq k \leq n.$

Consider the following sequence $\{g_n(x)\}_{n\geq 0}$ of functions 
such that $g_n:\BR^+ \to \BR^+$ for every $n.$
Let $g_0(x)\equiv 1,$ and define $g_1(x),g_2(x),\ldots$
inductively by letting
\begin{equation} \label{eqn:indg}
g_{n+1}(x)
=\int_{0}^{x}g_n(y)dy+\int_{x}^\infty e^{x-y}g_n(y)dy,
\end{equation}
for every $n\geq 1.$
\begin{prop} \label{prop:fcniter}
With definitions as above, we have that
\[
g_n(x)=\sum_{k=0}^n c_{k,n}\frac{x^k}{k!}.
\]
\end{prop}
\noindent
{\bf Proof.}
We start by noting that
\[
g_1(x)=\int_0^x 1 dy+e^x\int_{x}^\infty e^{-y}dy=x+1,
\]
and since $c_{1,1}=c_{0,1}=1$ the statement holds for $n=1.$
Assume therefore that it holds for $n-1$ and observe that with 
$c_{-1,n}=0$,
\begin{eqnarray*}
\lefteqn{g_n(x)=\int_0^x g_{n-1}(y) dy+e^x\int_{x}^\infty g_{n-1}(y)e^{-y}dy}\\
& & =\sum_{k=0}^{n-1}c_{k,n-1}\left(\int_0^x \frac{y^{k}}{k!} dy+e^x\int_{x}^\infty \frac{y^{k}}{k!}e^{-y}dy\right)\\
& & =\sum_{k=0}^{n-1}c_{k,n-1}\left(\frac{x^{k+1}}{(k+1)!}+\frac{x^{k}}{k!}+\cdots+1\right)\\
& & =\sum_{k=0}^{n-1}c_{k,n-1}\sum_{l=0}^{k+1}\frac{x^l}{l!}
=\sum_{k=0}^{n}\frac{x^{k}}{k!}\sum_{l=k-1}^{n-1}c_{l,n-1}.
\end{eqnarray*}
By using \eqref{eqn:cattri}, we conclude the proof.
\fbox{}\\

Our next result provides a link between the particle process $\zeta$ 
defined in the 
introduction, and the functions $g_n(x).$ Recall the interpretation that a particle
at position $Z_{k,n}=x,$ independently gives rise to new particles according to a Poisson process 
with intensity measure $d\mu_x=I(y\geq 0)ue^{-(y-x)^+}dy$, 
so that in particular the entire process is restricted to $\BR^+$. 
 Recall also the definition of 
$X^n_{[a,b]}$ in \eqref{eqn:Xndef}.

\begin{prop} \label{prop:contbranch}
Let
\[
F_n(R)=\BE[X^n_{[0,R]}].
\]
For any $u<\infty,$ $F_n(R)$ is differentiable with respect to $R$, 
and we have that with $f_n(R):=F'_n(R),$ 
\[
f_{n}(R)=u^{n}\sum_{k=0}^{n-1} c_{k,n-1}\frac{R^k}{k!}=u^n g_{n-1}(R),
\]
for every $n\geq 1.$
\end{prop}
\noindent
{\bf Proof of Proposition \ref{prop:contbranch}.}
We will 
prove the 
statement by induction, and so we start by noting that 
\[
F_1(R)=\BE[X^1_{[0,R]}]=uR,
\]
which follows since $Z_{1,0}$ is of type 0.
Therefore,  the statement holds for $n=1$.

Assume now that the statement holds for some fixed $n\geq 1.$
Let $R,\Delta R>0$ and consider 
\[
F_{n+1}(R+\Delta R)-F_{n+1}(R)=\BE[X^{n+1}_{[R,R+\Delta R]}].
\]
Any particle in generation $n$ of type smaller than $R$ gives rise to individuals
in $[R,R+\Delta R]$ (in generation $n+1$) at rate $u.$ Furthermore, any individual
of type $x\in [R,R+\Delta R]$ gives rise to individuals in $[R,R+\Delta R]$ at 
most at rate $u$ while individuals of type $x>R+\Delta R$ produce individuals
in $[R,R+\Delta R]$ at rate at most $ue^{R+\Delta R-x}.$ We therefore get the 
following upper bound
\begin{equation} \label{eqn:indUB}
\BE[X^{n+1}_{[R,R+\Delta R]}]
\leq u \Delta R \left(\BE[X^{n}_{[0,R+\Delta R]}]
+\sum_{k=0}^\infty\BE[X^{n}_{[R+\Delta R+k/N,R+\Delta R+(k+1)/N]}]e^{-k/N}\right),
\end{equation}
where $N$ is an arbitrary number. By assumption, $F_n(R)$ is differentiable,  
and by the mean value theorem,
\[
\BE[X^{n}_{[R+\Delta R+k/N,R+\Delta R+(k+1)/N]}]
\leq \frac{f_n(R+\Delta R+(k+1)/N)}{N},
\]
since $f_n(x)$ is increasing. Thus, we conclude from \eqref{eqn:indUB} that
\begin{eqnarray*}
\lefteqn{\BE[X^{n+1}_{[R,R+\Delta R]}]}\\
& & \leq \limsup_{N \to \infty} u \Delta R \left(F_{n}(R+\Delta R)
+\sum_{k=0}^\infty\frac{f_n(R+\Delta R+(k+1)/N)}{N} e^{-k/N}\right)\\
& & \leq \limsup_{N \to \infty}u \Delta R \left(F_{n}(R+\Delta R)
+\int_0^\infty\frac{f_n(R+\Delta R+(y+1)/N)}{N}e^{-(y-1)/N} dy\right)\\
& & =\limsup_{N \to \infty} u \Delta R \left(F_{n}(R+\Delta R)
+e^{1/N}\int_0^\infty f_n(R+\Delta R+z+1/N)e^{-z} dz\right)\\
& & =u \Delta R \left(F_{n}(R+\Delta R)
+\int_0^\infty f_n(R+\Delta R+z)e^{-z} dz\right),
\end{eqnarray*}
by the dominated convergence theorem. Hence, we conclude that 
\begin{eqnarray} \label{eqn:limsupFn}
\lefteqn{\limsup_{\Delta R \to 0}\frac{F_{n+1}(R+\Delta R)-F_{n+1}(R)}{\Delta R}}\\
& & \leq u\left(F_n(R)+\int_0^\infty f_n(R+z)e^{-z} dz\right)
=u\left(\int_0^Rf_n(z)dz+\int_R^\infty f_n(z)e^{R-z} dz\right),
\nonumber
\end{eqnarray}
again by the dominated convergence theorem.

Similarly, we get the following lower bound
\begin{eqnarray}\label{eqn:indLB}
\lefteqn{\BE[X^{n+1}_{[R,R+\Delta R]}]}\\
& & \geq u \Delta R \liminf_{N \to \infty} \left(\BE[X^{n}_{[0,R]}]
+\sum_{k=0}^\infty\BE[X^{n}_{[R+k/N,R+(k+1)/N]}]e^{-(k+1)/N}\right) \nonumber \\
& & \geq u \Delta R \liminf_{N \to \infty} \left(F_n(R)
+\sum_{k=0}^\infty\frac{f_n(R+k/N)}{N}e^{-(k+1)/N}\right) \nonumber \\
& & \geq u \Delta R \liminf_{N \to \infty} \left(F_n(R)
+\int_0^\infty\frac{f_n(R+(y-1)/N)}{N}e^{-(y+1)/N}dy\right) \nonumber \\
& & = u \Delta R \left(F_n(R)
+\int_0^\infty f_n(R+z)e^{-z}dz\right) \nonumber,
\end{eqnarray} 
which together with \eqref{eqn:limsupFn} gives us 
\[
\lim_{\Delta R \to 0}\frac{F_{n+1}(R+\Delta R)-F_{n+1}(R)}{\Delta R}
=u\left(\int_0^Rf_n(z)dz+\int_R^\infty f_n(z)e^{R-z} dz\right).
\]
Thus, we conclude that $F_{n+1}(R)$ is differentiable and that 
\[
f_{n+1}(R)=u^{n+1}\left(\int_0^Rg_{n-1}(z)dz+\int_R^\infty g_{n-1}(z)e^{R-z} dz\right)
=u^{n+1}g_n(R),
\]
where the last equality follows from \eqref{eqn:indg}.
\fbox{}\\

\noindent
{\bf Remarks:} The proof shows that for $u=1$, the functions 
$f_n(x)=F'_n(x)$ satisfies \eqref{eqn:indg}, which is of course why 
\eqref{eqn:indg} is introduced in the first place.

For future reference, we observe that $F_n(R)$ in fact depends on $u,$
and we sometimes stress this by writing $F_n(R,u)$. Furthermore, it is 
easy to see that for any $0<u<\infty$, we have that 
$F_n(R,u)=u^nF_n(R,1)$ for every $n\geq 1.$

\medskip

We have the following result
\begin{prop}\label{p.fupperbound}
Let $u<1/4,$ then  for every $x\geq 0,$
\[
\sum_{n=1}^{\infty}f_n(x)\leq u\frac{e^{4ux}}{1-4u}.
\]
\end{prop}
\noindent
{\bf Proof.}
By Propositions \ref{prop:fcniter} and \ref{prop:contbranch},
\begin{eqnarray}\label{eqn:fxu}
\lefteqn{\sum_{n=1}^{\infty}f_n(x)
=\sum_{n=1}^{\infty}u^{n}g_{n-1}(x)
=\sum_{n=0}^{\infty}u^{n+1}g_{n}(x)}\\
& & =\sum_{n=0}^{\infty}u^{n+1}\sum_{k=0}^n c_{k,n}\frac{x^k}{k!}
=\sum_{k=0}^{\infty}\frac{x^k}{k!}\sum_{n=k}^\infty u^{n+1}c_{k,n}. \nonumber
\end{eqnarray}
Furthermore, by using that ${m \choose n}$ is increasing in $m\geq n,$
we see that 
\begin{equation}\label{eqn:claim}
c_{k,n}=\frac{k+1}{n+1}{2n-k \choose n}
\leq {2n \choose n}
\leq \sum_{l=0}^{2n} {2n \choose l}=4^n.
\end{equation}
Combining \eqref{eqn:fxu} and \eqref{eqn:claim}, we see that for 
$u<1/4,$
\begin{equation}\label{e.fupperbound}
\sum_{n=1}^{\infty}f_n(x)\leq
\sum_{k=0}^{\infty}\frac{x^k}{k!}u\sum_{n=k}^\infty (4u)^n
=\frac{u}{1-4u}\sum_{k=0}^{\infty}\frac{(4ux)^k}{k!}=u\frac{e^{4ux}}{1-4u},
\end{equation}
finishing the proof. \fbox{} \\

\noindent
{\bf Remark:} As pointed out to us by an anonymous referee, a variant 
of Proposition \ref{p.fupperbound} can be proved along the following lines. 
Let $T$ be the integral operator defined by 
\[
T(g)=\int_0^x g(y)dy +\int_x^\infty e^{x-y}g(y)dy. 
\]
It is easy to check that $g(x)=(x+2)e^{x/2}$ is an eigenfunction of $T$
satisfying $T(g)=4g.$ Thus, since $g_0(x)\equiv 1\leq g(x)$ we get that
$g_1=T(g_0)\leq T(g)=4g,$ and iterating we see that 
$g_{n+1}=T(g_n)\leq T(4^n g)=4^{n+1}g.$ This can then be used in conjunction
with Proposition \ref{prop:contbranch} to prove the desired result. 

The justification for obtaining and using the explicit forms of 
$g_n, f_n$ and $F_n$ is twofold. Firstly, these forms will be 
convenient when proving the second part of Theorem \ref{thm:bp} and 
also when proving Lemma \ref{lemma:mballest} below. Secondly, 
we believe that the infinite type branching process $\zeta$ is of 
independent interest, and therefore a detailed analysis is 
intrinsically of value.

\medskip

\noindent
We can now prove Theorem \ref{thm:bp}.\\
\noindent
{\bf Proof of Theorem \ref{thm:bp}.}
We have that 
\[
F_n(R)=\int_0^Rf_n(x)dx,
\]
so that 
\[
\sum_{n=1}^\infty F_n(R)=\sum_{n=1}^\infty\int_0^Rf_n(x)dx.
\]
Furthermore, for $u<1/4$, we can use Proposition \ref{p.fupperbound} and
the dominated convergence theorem to conclude that 
\[
\sum_{n=1}^\infty F_n(R)\leq \int_0^R u\frac{e^{4ux}}{1-4u} dx
\leq\frac{e^{4uR}}{4(1-4u)}.
\]
We can now use Propositions \ref{prop:fcniter} and 
\ref{prop:contbranch} to get that 
\begin{eqnarray*}
\lefteqn{F_{n+1}(R)=u^{n+1}\int_0^Rg_n(x)dx
=u^{n+1}\int_0^R \sum_{k=0}^n c_{k,n}\frac{x^k}{k!}dx}\\
& & =u^{n+1}\sum_{k=0}^n c_{k,n}\frac{R^{k+1}}{(k+1)!}
\geq u^{n+1}c_{0,n}R
=\frac{u^{n+1}}{n+1}{2n \choose n}R
\geq \frac{u^{n+1}4^n}{2(n+1)^2}R,
\end{eqnarray*}
by using that $2(n+1){2n \choose n}\geq \sum_{l=0}^{2n}{2n \choose l}=4^n$
which follows since $l=n$ maximizes ${2n \choose l}.$ We see that if 
$u>1/4$, the right hand side diverges, and 
so the statement follows.
\fbox{}\\

\section{Proof of $u_c(d)>0.$} \label{sec:lb}
The aim of this section is to prove the lower bound of Theorem \ref{thm:main}. 
We will do this by establishing a link between the cylinder process $\omega$
and the particle process of Section \ref{s.actbt}. As an intermediate step, we 
will in Section \ref{subsec:mbp} consider particle processes with 
offspring distributions that can be weakly bounded above by $\zeta.$ 
In Section \ref{subsec:icp}, these new particle processes and
the cylinder process in $\BH^d$ will be compared.
 Thereafter, this link
is used in Section \ref{subsec:lb} to obtain the required lower bound.

\subsection{Particle processes weakly dominated by $\zeta$} \label{subsec:mbp}
Recall that $d\mu_{x}(y)=1_{(y\ge 0)}u e^{-(x-y)^+}\,dy$ and
suppose that $(\nu_{x})_{x\in {\mathbb R}^{+}}$ is a family of measures 
with the following property: there is a constant $c\in (0,\infty)$ such 
that for all integers $k,l\ge 0$
\begin{equation}\label{e.compeq}
\sup_{x\in(l,l+1]}\nu_x((k,k+1])\le c\inf_{x\in(l,l+1]}\mu_x((k,k+1]), 
\end{equation}
and moreover, $\nu_x(\{0\})=0$ for all $x\ge 0$. 
This last assumption is made only for convenience; if one allows the 
measures to have an atom at $0$ what follows below can be modified 
fairly easy to get similar conclusions.
The particle processes that we consider here are defined as the one in
Theorem \ref{thm:bp}, but using
$\nu_x$ as the offspring distribution in place of $\mu_x$ for a particle
of type $x.$ Recall that we think of the position of a particle 
in $\BR^+$ as being the type of that particle.
Of course, we still assume that every particle produces offspring independently. 
For this process, let $\tilde{X}_{D}^{n}$ be the number of individuals 
in generation $n$ of type in $D\subset {\mathbb R}_+$. Furthermore let
\[
\tilde{F}_n(R)={\mathbb E}[\tilde{X}_{[0,R]}^{n}].
\]
\begin{lem}\label{l.altproc}
With $c\in(0,\infty)$ as in \eqref{e.compeq}, we have that for every $R\in {\mathbb N}_+$
$$\tilde{F}_n(R)\le c^{n} F_n(R).$$
\end{lem}

{\bf Proof.}
Let $R\in {\mathbb N}_+$. It suffices to show that with 
$c$ as in \eqref{e.compeq}, and any integers $n\ge 1$ and $k\ge 0$, we have 
\begin{equation}\label{e.nts}
\BE[\tilde{X}^n_{(k,k+1]}]\le c^n \,\BE[X^n_{(k,k+1]}].
\end{equation}
Since $\tilde{F}_n(0)=F_n(0)=0$, the claim of the lemma will then follow by summing the two sides 
of~\eqref{e.nts} from $k=0$ to $k=R-1$.
We proceed by induction in $n$. For any $k\in {\mathbb N}$ we have that 
\[
\BE[\tilde{X}^1_{(k,k+1]}] 
=\nu_0((k,k+1])\stackrel{\eqref{e.compeq}}{\le} c\, \mu_0((k,k+1]) 
=c\,\BE[X^1_{(k,k+1]}],
\]
so that \eqref{e.nts} holds for $n=1.$
Assume therefore that \eqref{e.nts} holds for some $n\geq 1$ and every $k\geq 0.$
Let $\tilde{Y}_{k,l}^{n}$ denote the number of individuals in 
generation $n$ of type in $(k,k+1]$ with parents of type in $(l,l+1]$. We have
\begin{eqnarray}
\lefteqn{\BE[\tilde{X}^{n+1}_{(k,k+1]}]} \\
& & =\sum_{l=0}^\infty \BE[\tilde{Y}^{n+1}_{k,l}]
\le\sum_{l=0}^\infty \sup_{x\in(l,l+1]} \nu_x((k,k+1])\BE[\tilde{X}^n_{(l,l+1]}]\nonumber\\ & &
\stackrel{~\eqref{e.nts}}{\le} \sum_{l=0}^{\infty}c\, \inf_{x\in (l,l+1]}\mu_x((k,k+1])
c^{n} \BE[X^n_{(l,l+1]}]\le \,c^{n+1} \BE[X^{n+1}_{(k,k+1]}].\nonumber
\end{eqnarray}
This finishes the proof of the lemma.
\fbox{}\\

\subsection{The independent cylinder process} \label{subsec:icp}
We now turn to the independent cylinder process discussed in the introduction. 
We start by defining the process itself, and the coupling with the 
ordinary line process
$\omega.$ Thereafter we establish a link between our independent 
cylinder process and 
the particle process studied in Section \ref{subsec:mbp}.

Formally, we define the independent cylinder process as follows. 
Let $\omega, (\omega_{k,n})_{k,n=1}^\infty$ be an i.i.d.
collection of Poisson line processes with intensity $u\mu_{d,1}.$ 
We use $\omega$ to define 
$\CC(\omega)$. Fix any (deterministic) line 
$L_{1,0}$ such that $o\in L_{1,0}.$ This is the single line of generation $0.$
Let 
\[
\eta^{1}:=\{L\in \omega:\cc(L)\cap \cc(L_{1,0})\neq \emptyset\}.
\]
Thus $\eta^1$ is simply the collection of lines in $\omega$ 
such that the corresponding 
cylinders intersect $\cc(L_{1,0}).$ 
Recall the definition of $\rho(L)$ for $L\in A(d,1)$ from 
Section \ref{s.geodesics}.
Let $L_{1,1},L_{2,1},\ldots$ be 
the enumeration of the lines in $\eta^1$ satisfying
$\rho(L_{k,1})<\rho(L_{k+1,1})$ for every $k\geq 1.$ As when we defined 
the particle process $\zeta$ in the introduction, the particular choice 
of enumeration is somewhat arbitrary. Define
\[
\eta^{1,2}:=\{L\in\omega\setminus\eta^1: 
\cc(L)\cap \cc(L_{1,1})\neq \emptyset\}\cup 
\{L\in \omega_{1,1}:\cc(L)\cap \cc(L_{1,1})\neq \emptyset, 
\cc(L)\cap \cc(L_{1,0})\neq \emptyset\}.
\]
The first set of lines corresponds to the cylinders in $\omega$ that intersect
$\cc(L_{1,1})$ but are not included in the definition of $\eta^1$ (i.e. intersect
$\cc(L_{1,0})).$  The second
set of lines is an independent copy of the set of lines in $\omega$ that 
intersects both $\cc(L_{1,0})$ and $\cc(L_{1,1}).$ Thus, we see that $\eta^{1,2}$
and $\eta^{1}$ are created in the same way, i.e. by considering the set of cylinders
of a Poisson cylinder process intersecting $\cc(L_{1,1})$ and $\cc(L_{1,0})$
respectively. For any $k\geq 1,$ let 
\begin{eqnarray*}
\lefteqn{\eta^{k,2}:=\{L\in\omega\setminus\left(\eta^1\cup_{l=1}^{k-1}\eta^{l,1}\right): 
\cc(L)\cap \cc(L_{k,1})\neq \emptyset\}}\\
& & \ \ \ \ \ \ \cup 
\{L\in \omega_{k,1}:\cc(L)\cap \cc(L_{k,1})\neq \emptyset, 
\cc(L)\cap \left(\cc(L_{1,0})\cup_{l=1}^{k-1}\cc(L_{l,1})\right)\neq \emptyset\}.
\end{eqnarray*}
We think of $\eta^{k,2}$ as being created from $\omega$ where $\omega$
has not already been used, and from  $\omega_{k,1}$
where $\omega$ has been used. From this construction it is obvious that given $\eta^1,$
the sequence $\eta^{1,2}, \eta^{2,2},\ldots $ is independent.
We let 
\[
\eta^2=\cup_{k=1}^\infty \eta^{k,2},
\]
and let $L_{1,2},L_{2,2},\ldots$ be the enumeration of the lines in 
$\eta^2$ satisfying $\rho(L_{1,2})<\rho(L_{2,2})<\cdots$. 
These are the lines belonging to generation 2.

We proceed in the obvious
way, defining $\eta^{k,n}$ and $\eta^n$ for $k,n\geq 1.$ Finally, let 
\[
\eta:=\cup_{n=1}^\infty \eta^n,
\]
and
\begin{equation} \label{eqn:defC0eta}
\CC_0(\eta)=\cc(L_{1,0})\bigcup_{L\in \eta}\cc(L).
\end{equation}
Consider the set $\cc(L_{1,0})\cup \CC(\omega),$ and 
define $\CC_0(\omega) \subset \cc(L_{1,0})\cup \CC(\omega)$ to be the 
maximally connected component of the origin.
By our construction, 
\begin{equation}\label{eqn:C0coupling}
\CC_0(\omega) \subset \CC_0(\eta).
\end{equation}Indeed, any cylinder 
in $\CC_0(\omega)$ intersecting $L_{1,0}$ is by definition in $\eta^1,$ and 
in general, any cylinder in $\CC_0(\omega)$ separated from $L_{1,0}$
by $k$ other cylinders will belong to $\eta^{k+1}.$

From $\eta$, a particle process $\bar{\zeta}$ is induced in the following way.
With $L_{1,n},L_{2,n},\ldots$ being the enumeration of $\eta^n$ 
satisfying $\rho(L_{1,n})<\rho(L_{2,n})<\cdots,$ we let 
$\bar Z_{k,n}=\bar Z_{k,n}(L_{k,n}):=\rho(L_{k,n})$ 
for every $k=1,2,\ldots.$ Furthermore, 
$\bar \zeta^{k,n}=\{\bar Z_{l,n}(L_{l,n}): L_{l,n}\in \eta^{k,n}\},$ and 
of course $\bar \zeta^n=\bigcup_{k=1}^\infty \bar \zeta^{k,n}.$
Since given $\eta^n,$ the sequence $\eta^{1,n+1},\eta^{2,n+1},\ldots$ 
is independent, 
it follows that $\bar \zeta^{1,n+1}, \bar \zeta^{2,n+1},\ldots$ are 
independent given 
$\bar \zeta^n.$ Therefore, $\bar \zeta=(\bar \zeta^n)_{n=1}^\infty$ 
has the desired independence 
properties. We note the similarities between $\zeta$ and $\bar \zeta$,
and that the only essential difference lies in the offspring distributions,
which we address next.

Fix some $L_{k,n}\in \eta^n$ and consider an offspring $L\in \eta^{k,n+1}.$ 
If $L$ is
at distance between $l$ and $l+1$ from the origin, then this corresponds to 
an offspring $\bar Z\in \bar \zeta^{k,n+1}$ of $\bar Z_{k,n}$ such 
that $\bar Z\in (l, l+1].$ 
Furthermore, the expected number of offspring (of $L_{k,n}$) belonging to 
$\CL_{B(o,l+1)}\setminus \CL_{B(o,l)}$ equals 
$u\mu_{d,1}\left( \CL_{\cc(L_{k,n},2)}\cap 
(\CL_{B(o,l+1)}\setminus \CL_{B(o,l)}) \right),$
and so we see that the particle process $\bar \zeta$ 
can be described using the intensity measures
$\{\tau_x\}_{x\geq 0}$ where
\begin{equation}\label{e.taufact}
\tau_x((l,l+1])
=u \mu_{d,1}\left( \CL_{{\mathfrak c}(L_x,2)}
\cap (\CL_{B(o,l+1)}\setminus \CL_{B(o,l)}) \right),
\end{equation}
with $L_x$ satisfying $x=d_H(o,L_x).$ Our next result will be used to 
prove that $\{\tau_x\}_{x\geq 0}$ satisfies \eqref{e.compeq} for some
$c<\infty.$

\begin{lem}\label{l.skallemma}
Let $x\in \BR^+$ and $L\in A(d,1)$ be such that $d_H(o,L)=x.$ There exists a constant
$C(d)\in(0,\infty)$ such that for any $k\geq 0,$
\[
\mu_{d,1}(\CL_{{\mathfrak c}(L,2)}\cap \left(\CL_{B(o,k+1)}\setminus \CL_{B(o,k)})\right)
\le C \exp(-(d-1)(x-k)^+).
\]
\end{lem}

{\bf Proof.}
Fix $k\in {\mathbb N}$. Suppose that $L=\{\gamma(t)\,:\,-\infty<t<\infty\}$ where 
the parametrization of $\gamma$ is chosen to be unit speed and so that 
$d_H(o,L)=d_H(o,\gamma(0))$. For $i\in{\mathbb Z}$ let $y_i=\gamma(i)$ 
and $B_i=B(y_i,3)$. Then ${\mathfrak c}(L,2)\subset \cup B_i$ since any point
in $\cc(L,2)$ is at distance at most $2$ from $L,$ and any point in $L$ is at distance
at most $1/2$ from some $y_i$. We now claim that 
\begin{equation}\label{e.ybound}
d_H(o,y_{i+1})-d_H(o,y_i)\ge c_1
\end{equation} 
for every $i\geq 0,$ and some constant $c_1>0$. To see this, assume that 
$i \geq 0,$ 
and observe that by \eqref{eqn:CS1},
\[
d_H(0,y_i)=\cosh^{-1}(\cosh(x)\cosh(i)),
\]
since the angle between $L$ and the geodesic from $o$ to $L$ is $\pi/2.$
Equivalently, we get that 
\[
d_H(0,y_i)=\log\left(\cosh(x)\cosh(i)+\sqrt{\cosh^2(x)\cosh^2(i)-1}\right).
\]
Hence, 
\begin{eqnarray*}
\lefteqn{d_H(0,y_{i+1})-d_H(0,y_i)}\\ 
& &=\log\left( \frac{\cosh(x)\cosh(i+1)+\sqrt{\cosh^2(x)\cosh^2(i+1)-1}}
{\cosh(x)\cosh(i)+\sqrt{\cosh^2(x)\cosh^2(i)-1}}\right)\\ 
& & \ge \log\left( \frac{\cosh(x)\cosh(i)\cosh(1)+\sqrt{\cosh^2(x)\cosh^2(i)\cosh^2(1)-1}}
{\cosh(x)\cosh(i)+\sqrt{\cosh^2(x)\cosh^2(i)-1}}\right)\\ 
& & \ge \log\left( \frac{\cosh(x)\cosh(i)\cosh(1)+\sqrt{\cosh^2(x)\cosh^2(i)\cosh^2(1)-\cosh^2(1)}}
{\cosh(x)\cosh(i)+\sqrt{\cosh^2(x)\cosh^2(i)-1}}\right)\\ 
& & =\log(\cosh(1)),
\end{eqnarray*}
where we use that $\cosh(i+1)\geq \cosh(i)\cosh(1)$ which holds since $i\geq 0.$
Hence, ~\eqref{e.ybound} follows with $c_1=\log(\cosh(1))$.

Assume first that $k<x$. From \eqref{e.ybound}, we
get that $d(y_i,0)\ge x+c_1 |i|$ for every $i$ using symmetry. We get that
\begin{eqnarray*}
\lefteqn{\mu_{d,1}(\CL_{{\mathfrak c}(L,2)}\cap \left(\CL_{B(o,k+1)}\setminus \CL_{B(o,k)})\right)}\\
& & \le \mu_{d,1}(\CL_{{\mathfrak c}(L,2)}\cap \left(\CL_{B(o,k+1)})\right) 
\le \sum_{i\in {\mathbb Z}}\mu_{d,1}(\CL_{B_i}\cap \CL_{B(o,k+1)})\\
& & \le \sum_{i\in {\mathbb Z}} C \exp(-(d-1) (x+c_1 |i|-k)
  \le C' \exp(-(d-1)(x-k)),
\end{eqnarray*}
where the penultimate inequality follows from Lemma~\ref{l.diffballs}. 
Now assume instead that $x\le k$.  Let $p=\inf\{|i|\,:\,d_H(o,y_i)\ge k-3\}$. 
Using the union bound and that  
$\mu_{d,1}(\CL_{B_i}\cap \left(\CL_{B(o,k+1)}\setminus \CL_{B(o,k)})\right)=0$ 
when $|i| \le p$, we get
\begin{eqnarray*}
\lefteqn{\mu_{d,1}(\CL_{{\mathfrak c}(L,2)}\cap \left(\CL_{B(o,k+1)}\setminus \CL_{B(o,k)})\right)}\\
& & \le \sum_{i\in {\mathbb Z}}\mu_{d,1}(\CL_{B_i}\cap \left(\CL_{B(o,k+1)}\setminus \CL_{B(o,k)})\right)\\
& & \le 2\sum_{ i=p}^\infty \mu_{d,1}(\CL_{B_i}\cap \CL_{B(o,k+1)})\le 2 \sum_{i=p}^\infty c \exp(-(d-1)(c_1(i-p))) \le C,
\end{eqnarray*}
where again we use Lemma~\ref{l.diffballs}.
\fbox{}\\

Let
\[
\bar Z_{[0,R]}^n=\sum_{k=1}^\infty I(\bar Z_{k,n}\leq R),
\]
and let $H_n(R,u)=\BE[\bar Z^n_{[0,R]}].$ It is not hard to show that, 
similarly to the observation after the proof of Proposition 
\ref{prop:contbranch}, $H_n(R,u)=u^n H_n(R,1).$ This follows from the
Poissonian nature of the process.

We can now use Lemma \ref{l.skallemma}
to show the following result
(recall the definition of $F_n(R)=F_n(R,u)$ from Section~\ref{s.actbt}).
\begin{lem}\label{l.FHcompare}
There is a constant $c(d)\in (0,\infty)$ such that for every 
$R\in {\mathbb N}_{+}$, and $0<u<\infty,$
$$H_n(R,u)\le c(d)^n F_n(R,u).$$
\end{lem}
{\bf Proof.}
Since $H_n(R,u)=\BE[\bar Z^n_{[0,R]}]$, and the particle process 
$\bar \zeta$
uses $\{\tau_x\}_{x\geq 0}$ as intensity measures, it suffices in view
of Lemma~\ref{l.altproc}, 
to show that there is a constant $c<\infty$ such that for every integer $k,l\ge 0$, 
\begin{equation} \label{eqn:msrineq}
\sup_{x\in(k,k+1]}\tau_x((l,l+1])\le c \inf_{x\in (k,k+1]}\mu_x((l,l+1]).
\end{equation}
From Lemma~\ref{l.skallemma}, we have 
\begin{eqnarray} \label{e.tau1}
\lefteqn{ \sup_{x\in(k,k+1]}\tau_x((l,l+1])}\\
& & \le c u \sup_{x\in(k,k+1]} \exp(-(d-1)(x-l)^+) \le c' u\exp(-(d-1)(k-l)^+).\nonumber
\end{eqnarray}
On the other hand
\begin{equation}\label{e.mu1}
\inf_{x\in(k,k+1]}\mu_x((l,l+1])=\inf_{x\in(k,k+1]}u\int_{l}^{l+1} \exp(-(x-y)^+)\,dy\ge c u\exp(-(k-l)^+).
\end{equation}
Equations~\eqref{e.tau1} and~\eqref{e.mu1} establishes \eqref{eqn:msrineq},
and the lemma follows.
\fbox{}\\
In what follows, we drop the explicit dependence on $u$ from the notation
and simply write $H_n(R)$ and $F_n(R).$

\subsection{Proof of Theorem \ref{thm:main}} \label{subsec:lb}
We now have all the ingredients to prove our main result.

\noindent
{\bf Proof of Theorem \ref{thm:main}.}\\
Using Lemma \ref{lemma:mon} and Proposition \ref{p.cylconn}, we only need to 
prove that $u_c(d)>0$. To that end, let 
\[
V(R)=\BE\left[|\{L\in \CC_0(\omega):L\cap B(o,R)\neq \emptyset\}|\right],
\] 
that is, $V(R)$ is the expected number of cylinders in $\CC_0(\omega)$
which intersect $B(o,R)$. Recall that $\CC_0(\omega)$ is the maximally connected
component of $\cc(L_{1,0})\cup \CC(\omega),$ and recall also the definition of 
$\CC_0(\eta)$ from \eqref{eqn:defC0eta}. By \eqref{eqn:C0coupling} 
we can couple $\CC_0(\omega)$ and $\CC_0(\eta)$
so that $\CC_0(\omega)\subset \CC_0(\eta),$ and so we have as in the proof of Theorem 
\ref{thm:bp}, that for $u<1/(4c)$ with $c=c(d)$ as in Lemma ~\ref{l.FHcompare},
\[
V(R)\le \sum_{n=0}^{\infty}H_n(R)\leq \sum_{n=0}^{\infty} c^n F_n(R)
\leq \frac{Ce^{4ucR}}{4(1-4cu)}
\]
by using Lemma~\ref{l.FHcompare} in the second inequality.

Hence, we see that when $0<u<1/(4c)$, $V(R)$ grows at most exponentially in $R$ 
at a rate which is strictly smaller than $1$. On the other hand, we have that 
\[
\BE[|L\in \omega:L\cap B(o,R)|]=u\mu_{d,1}(\CL_{B(o,R)})
=Cu\sinh^{d-1}(R),
\]
by \eqref{eqn:linevolume}. This grows exponentially at rate $(d-1)R$ and so we see that
with probability one,  
$\CC_0(\omega)$ is a strict subset of $\cc(L_{1,0})\cup \CC(\omega).$ 
We conclude that $\CC(\omega)$ is a.s. not connected for this choice of $u.$
\fbox{}\\

\section{Proof of Theorem \ref{thm:infdiam}.} \label{sec:infdiam}
Similar to the notation of Section \ref{sec:prel} we let 
$A \stackrel{m}{\leftrightarrow} B$ denote the event 
that there exists $1\leq l\leq m$ and a sequence
of cylinders $\cc^1,\cdots,\cc^l\in\omega$ such that 
$A\cap \cc^1 \neq \emptyset, \cc^1\cap \cc^2 \neq \emptyset,
\ldots, \cc^l\cap B\neq \emptyset.$
That is, the sequence $\cc^1,\cdots,\cc^l$ connects $A$ to $B$ in  $l$ steps.
We observe that $\{A \leftrightarrow B\}
=\{A \stackrel{1}{\leftrightarrow} B\}$.

We start with the following lemma.
\begin{lem} \label{lemma:mballest}
There exists a constant $D(m)<\infty$ (depending only on $u$ and 
$m$) such that for any $x,y\in \BH^d$ we have that 
\begin{equation} \label{eqn:mballest}
\BP[B(x,1) \stackrel{m}{\leftrightarrow} B(y,1)]
\leq D(m)d_H(x,y)^me^{-(d-1)d_H(x,y)}.
\end{equation}
\end{lem}
{\bf Proof.}
Assume without loss of generality that $x=o$ and let $d_H(o,y)=R.$
As in Section \ref{subsec:icp}, the expected number 
of cylinders in $\eta$ 
up to generation $m$ that intersect the ball $B(o,R)$ is bounded by 
\begin{eqnarray*}
\lefteqn{\sum_{n=1}^m H_n(R)\leq \sum_{n=1}^m c(d)^nF_n(R)}\\
& & = \sum_{n=1}^m c(d)^n\int_0^R f_n(x)dx
=\sum_{n=1}^m c(d)^nu^n\int_0^R g_{n-1}(x)dx \\
& & =\sum_{n=1}^m c(d)^nu^n\int_0^R \sum_{k=0}^{n-1}c_{k,n-1}\frac{x^k}{k!}dx
 \leq \sum_{n=1}^m c(d)^nu^nC4^n\int_0^R \sum_{k=0}^{n-1}\frac{x^k}{k!}dx\\
& & =\sum_{n=1}^m c(d)^nu^nC4^n\sum_{k=1}^{n}\frac{R^k}{k!}
=\sum_{k=1}^{m}\frac{R^k}{k!}C\sum_{n=k}^m (4cu)^n
\leq D(m)R^m,
\end{eqnarray*}
for some constant $D(m)<\infty$.

There exists a collection $\CB_R$ of balls of radius $1/4$ with centers in $\partial B(o,R)$
such that $|\CB_R|\geq ce^{(d-1)R}$ for some $c>0$ and such that 
any cylinder intersecting $B(o,R)$ intersects at most $c(d)<\infty$ balls in $\CB_R.$ 
To construct such a collection $\CB_R$, we consider first $D$ as in Lemma \ref{lemma:D}.
Let $G_R=D\cap (B(o,R+3/2)\setminus B(o,(R-1/2)^+))$. By a slight modification 
of the lower bound in \eqref{e.duniform}, we get that 
$|G_R|\ge c e^{(d-1) R}/\nu_d(B(o,1/2))=c' e^{(d-1)R}$.
For $x\in G_R$, let $x'$ be defined as the point on $\partial B(o,R)$ such that $x'$ 
minimizes the distance to $\partial B(o,R)$, and let $G_R'$ be the collection of all such $x'$. 
Obviously, the collection of balls $\CB_R:=\{B(x,1/4)\}_{x\in G_R'}$ satisfies
$|\CB_R|\ge c' e^{(d-1) R}$. Now let $L$ be a line intersecting $B(o,R+5/4)$ 
(only cylinders centered around such lines might intersect some ball in $\CB_R$). 
Using \eqref{e.ybound}, there is a universal constant $c_2<\infty$ and two 
points $x_1,x_2\in \partial B(o,R)$ (these points depend on $L$) such that 
${\mathfrak c}(L)\cap (B(o,R+1/4)\setminus B(o,(R-1/4)^+))\subset B(x_1,c_2)\cup B(x_2,c_2)$. 
Hence the number of balls from $\CB_R$ intersecting ${\mathfrak c}(L)$ is bounded by the 
number of points in $D\cap (B(x_1,c_2+2)\cup B(x_2,c_2+2)).$ This in turn is bounded 
by some constant $c_3(d)<\infty$, by the upper bound of \eqref{e.duniform}. 
Hence, the existence of the  $\CB_R$ is verified.

Using $\CB_R$, we see that the probability that a fixed ball at distance $R$ from 
$o$ will be intersected by any cylinder in $\eta$ of generation less than or 
equal to $m$ is bounded by 
\[
\frac{c_3}{|\CB_R|}\sum_{n=1}^m H_n(R)
\leq \frac{D(m)R^m}{e^{(d-1)R}},
\] 
by possibly increasing the value of $D(m).$

The statement follows by using that $\CC_0(\omega)\subset \CC_0(\eta)$
and noting that any cylinder that 
intersects $B(o,1)$ must also intersect the cylinder $\cc(L_{1,0}).$
\fbox{}\\

{\bf Proof of Theorem \ref{thm:infdiam}.}
Let $m\in {\mathbb N}_+$ and fix $\epsilon\in(0,1)$. 
Using Lemma \ref{lemma:mballest}, we can choose $r=r(m,\epsilon)<\infty$ 
so large that the probability that any two fixed cylinders separated by distance
$r$ will be connected in at most $m$ steps is less than $\epsilon.$
Indeed, take $r$ so large that the probability that 
$B(o,1)$ and $B(y,1)$ (where $y\in \partial B(o,r)$)
are connected in at most $m+2$ steps is less than $\epsilon$. Consider then two cylinders 
$\cc_1,\cc_2$ separated by distance $r,$ and assume without loss
of generality that $\cc_1\cap B(o,1)\neq \emptyset$ and $\cc_2 \cap B(y,1)\neq \emptyset.$
Then, if the probability that $\cc_1,\cc_2$ are connected in at most $m$ steps is  larger 
than $\epsilon,$ this would lead to a contradiction.

For lines $L_1,L_2\in A(d,1)$, let $E_m(L_1,L_2)$ be the event that 
${\mathfrak c}(L_1)$ and ${\mathfrak c}(L_2)$ are connected in at most
$m$ steps. Define the event
\[
H=\left\{\sum_{(L_1,L_2)\in \omega_{\neq}^{2}}I(E_m(L_1,L_2)^c)\ge 1\right\},
\] 
where the union is over all $2$-tuples of distinct lines in 
$\omega$. In words, $H$ is the event that there is at least one 
pair of lines in $\omega$ whose corresponding cylinders are 
not connected in at most $m$ steps.
We now let ${\mathbb E}^{L_1,L_2}$ denote the expectation with 
respect to $\omega+\delta_{L_1}+\delta_{L_2}$. Using the 
Slivnyak-Mecke formula (\cite{SchWeil} Corollary $3.2.3$) 
we get that
\begin{eqnarray*}
\lefteqn{{\mathbb E}\left[\sum_{(L_1,L_2)\in \omega^2_{\neq}} I(E_m(L_1,L_2)^c)\right]}\\
& & =u^2\int_{A(d,1)} \int_{A(d,1)}{\mathbb E}^{L_1,L_2}\left[I(E_m(L_1,L_2)^c)\right] \mu_{d,1}(d L_1) \mu_{d,1}(d L_2) \\
& & = u^2\int_{A(d,1)} \int_{A(d,1)}{\mathbb P}^{L_1,L_2}\left[E_m(L_1,L_2)^c\right] \mu_{d,1}(d L_1) \mu_{d,1}(d L_2) \\
& &  = u^2\int_{A(d,1)} \int_{A(d,1)}{\mathbb P}\left[E_m(L_1,L_2)^c\right] \mu_{d,1}(d L_1) \mu_{d,1}(d L_2) \\
& & \ge u^2 \int_{A(d,1)} \int_{A(d,1)}I(d_H(L_1,L_2)\ge r){\mathbb P}\left[E_m(L_1,L_2)^c\right] \mu_{d,1}(d L_1) \mu_{d,1}(d L_2) \\
& & \ge u^2\int_{A(d,1)} \int_{A(d,1)} I(d_H(L_1,L_2)\ge r)(1-\epsilon)\mu_{d,1}(d L_1) \mu_{d,1}(d L_2).
\end{eqnarray*}
Obviously, the expression on the right hand side diverges for any $0<u<\infty,$
so that 
\[
{\mathbb E}\left[\sum_{(L_1,L_2)\in \omega^2_{\neq}} I(E_m(L_1,L_2)^c)\right]
=\infty,
\]
and from this, it follows that ${\mathbb P}[H]>0$. Since the event 
$H$ is invariant under isometries of ${\mathbb H}^d$, it follows from
Proposition~\ref{p.01law} that ${\mathbb P}[H]=1$. 
Hence for any $m\in {\mathbb N}_+$ we have 
${\mathbb P}[\textrm{diam}({\mathcal C})\ge m-2]=1$, from which 
we conclude that ${\mathbb P}[\textrm{diam}({\mathcal C})=\infty]=1$. 
\fbox{}\\

\section{Proof of Proposition \ref{prop:infnbrinfcomp}} \label{sec:inic}
In this section, we prove that when $u<u_c(d)$, there are a.s. infinitely many 
connected components in ${\mathcal C}.$ 
Let $N(\omega)=$the number of connected components in ${\mathcal C}$.

{\bf Proof of Proposition \ref{prop:infnbrinfcomp}}
Obviously, the event $N(\omega)=k$ is invariant under isometries of 
$\BH^d$ and so 
using  Proposition~\ref{p.01law}, we have that for any 
$u$ there is $k=k(u)\in {\mathbb N}\cup\{\infty\}$ such that 
${\mathbb P}[N(\omega)=k]=1$. Suppose $u<u_c(d)$ and suppose  
that $1<k(u)<\infty$.  It is not hard to show that there exist 
points $y_1,...,y_k \in {\mathbb H}^d$ such that the event 
\[
A:=\{\cup_{i=1}^k y_i\mbox{ intersects all components of }{\mathcal C}(\omega)\}
\cap \{\omega({\mathcal L}_{B(o,1)})=0\}
\] 
has positive probability. 

Now let $\omega_1$ be the restriction of $\omega$ to ${\mathcal L}_{B(o,1)},$ 
and let $\omega_2$ be the restriction of $\omega$ to $({\mathcal L}_{B(o,1)})^c$.
Since ${\mathbb P}[A]>0$, it follows that
\[
B:=\{\cup_{i=1}^k y_i\mbox{ intersects all components of }{\mathcal C}(\omega_2)\}
\] 
has positive probability. Define the event 
\[
C:=\{\cup_{i=1}^k y_i\subset {\mathcal C}(\omega_1)\}\cap\{{\mathcal C}(\omega_1)
\mbox{ is connected}\}.
\] 
It is easy to see that ${\mathbb P}[C]>0$ (indeed, $\omega_1$ might consist 
of $k$ lines $L_1$,...,$L_k$ such that ${\mathfrak c}(L_i)$ contains $o$ and $y_i$). 
Since $\omega_1$ and $\omega_2$ are independent, it follows that $B$ and $C$ are 
independent and hence ${\mathbb P}[B\cap C]>0$. The event $B\cap C$ 
implies that $N(\omega)=1$, whence 
$\BP[N(\omega)=1]>0$ which contradicts $\BP[N(\omega)=k]=1$.
We conclude that 
$N(\omega)\in \{1,\infty\},$ and since $u<u_c$ by assumption, it follows that a.s.
$N(\omega)=\infty$.
\fbox{}\\

{\bf Acknowledgements: }We would like to thank I. Benjamini for suggesting 
some of the problems dealt with in the paper. We would also like to thank the anonymous referees
for many useful suggestions and comments.

\end{document}